\documentclass[final]{siamart220329}   

\usepackage{siampaper}
\usepackage{comment}
\usepackage{caption,subcaption}
\usepackage[shortlabels]{enumitem}
\usepackage{booktabs}
\usepackage{multirow}
\usepackage{framed}
\usepackage{threeparttable}
\usepackage{xr}
\usepackage{bookmark}

\newcommand{\bmat}[1]{\begin{bmatrix} #1 \end{bmatrix}}

\newcommand*\circled[1]{\raisebox{.5pt}{\textcircled{\raisebox{-.9pt} {#1}}}}

\makeatletter
\def\thmt@innercounters{section,equation,theorem}
\makeatother

\makeatletter
\g@addto@macro\bfseries{\boldmath}
\makeatother

\usepackage[utf8]{inputenc}
\usepackage[T1]{fontenc}

\usepackage{tikzscale}
\usepackage{pgfplots}
\pgfplotsset{compat=newest}
\usetikzlibrary{backgrounds}
\usetikzlibrary{intersections}
\usetikzlibrary{external}
\usetikzlibrary{math}

\graphicspath{{graphics}}  

\makeatletter
\newcommand\fs@ruled@notop{
    \def\@fs@cfont{\bfseries} 
    \let\@fs@capt\floatc@ruled
    \def\@fs@pre{} 
    \def\@fs@post{} 
    \def\@fs@mid{} 
    \let\@fs@iftopcapt\iffalse
}
\renewcommand\fst@algorithm{\fs@ruled@notop} 
\makeatother

\makeatletter
\usepackage{etoolbox}
\patchcmd\H@refstepcounter{\protected@edef}{\protected@xdef}{}{}
\makeatother

\title{%
  Revisiting sparse matrix coloring and bicoloring
}

\author{%
  Alexis Montoison%
  \thanks{%
    Argonne National Laboratory, IL, USA.
    E-mail: \mailto{amontoison@anl.gov}
  }
  \and
  Guillaume Dalle%
  \thanks{%
    LVMT, ENPC, Institut Polytechnique de Paris, Univ Gustave Eiffel, Marne-la-Vallée, France.
    E-mail: \mailto{guillaume.dalle@enpc.fr}
  }
  \and
  Assefaw Gebremedhin%
  \thanks{%
    Washington State University, School of EECS, WA, USA.
    E-mail: \mailto{assefaw.gebremedhin@wsu.edu}\!\!\!\!
  }
  \indent Version of \today.%
}
\date{\today}

\headers{Revisiting sparse matrix coloring and bicoloring}{A. Montoison, G. Dalle and A. Gebremedhin}

\begin{document}

\maketitle

\thispagestyle{firstpage}
\pagestyle{myheadings}

\begin{abstract}
    Sparse matrix coloring and bicoloring are fundamental building blocks of sparse automatic differentiation. 
    Bicoloring is particularly advantageous for rectangular Jacobian matrices with at least one dense row and column.
    Indeed, in such cases, unidirectional row or column coloring demands a number of colors equal to the number of rows or columns.
    We introduce a new strategy for bicoloring that encompasses both direct and substitution-based decompression approaches.
    Our method reformulates the two variants of bicoloring as star and acyclic colorings of an augmented symmetric matrix.
    We extend the concept of neutral colors, previously exclusive to bicoloring, to symmetric colorings, and we propose a post-processing routine that neutralizes colors to further reduce the overall color count.
    We also present the Julia package \texttt{SparseMatrixColorings.jl}, which includes these new bicoloring algorithms alongside all standard coloring methods for sparse derivative matrix computation.
    Compared to \texttt{ColPack}, the Julia package also offers enhanced implementations for star and acyclic coloring, vertex ordering, as well as decompression.
\end{abstract}




\section{Introduction}

\subsection{Context and motivation}

Modern scientific computing makes heavy use of automatic differentiation (AD) to obtain derivatives of computer programs \citep{blondelElementsDifferentiableProgramming2024}.
Such derivatives are relevant in a wide variety of settings, from differential equations and nonlinear optimization to machine learning.
Whenever the objects of interest are sparse matrices (Jacobians or Hessians), their computation can be accelerated significantly \citep{griewankEvaluatingDerivativesPrinciples2008}.
However, reaping the benefits of sparsity comes with increased algorithmic and implementation challenges.

Consider a function $f : \mathbb{R}^n \to \mathbb{R}^m$ and its Jacobian matrix $J = \partial f(x)$ at a point $x$.
For roughly the same cost as the evaluation of $f$, forward-mode AD evaluates a single Jacobian-vector product $u \mapsto J u$, while reverse-mode AD evaluates a single Vector-Jacobian product $v \mapsto v^{\top\!} J$.
Importantly, this is obtained without ever materializing the matrix~$J$.
Taking the seed $u = e_j$ (resp. $v = e_i$) to be a basis vector of the input space (resp. output space) lets us recover the $j$-th column (resp. $i$-th row) of the Jacobian.
One can then iterate through basis vectors and reconstruct the full Jacobian, either column-wise or row-wise.
The complexity of this recovery, expressed in number of calls to $f$, scales either with the input dimension $n$ or with the output dimension $m$.

Luckily, when the Jacobian is sparse, focusing on the non-zero entries often speeds things up.
For instance, if the non-zeros in columns $j_1$ and $j_2$ never overlap at any row, they can be recovered together in a single product $J (e_{j_1} + e_{j_2})$ without the non-zero entries  getting mixed up.
Making this insight actionable requires solving a \emph{matrix
coloring problem}~\citep{gebremedhinWhatColorYour2005}, which determines the sets of columns or rows of $J$ that can be evaluated simultaneously.
The number of distinct colors obtained corresponds to the number of seeds required to recover all non-zero entries of $J$.
The recovery complexity then scales with the number of colors $c$, rather than $n$ or $m$.

\begin{figure}[ht]
    \centering
    \includegraphics[width=0.30\linewidth]{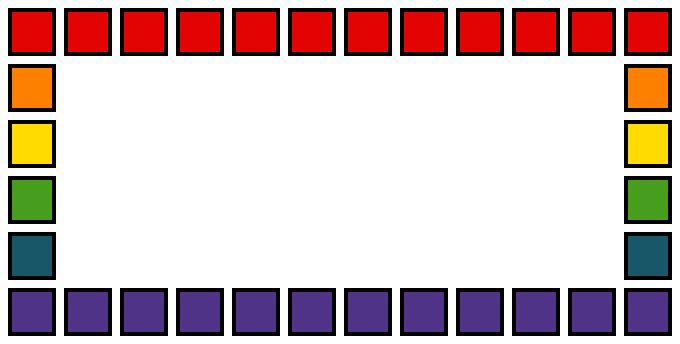}
    \hspace{2cm}
    \includegraphics[width=0.30\linewidth]{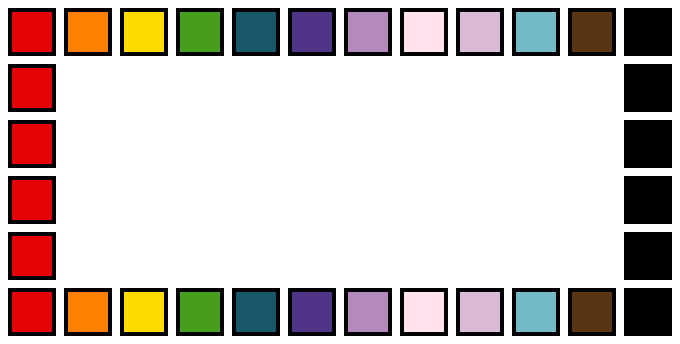}
    \caption{Row coloring (left) and column coloring (right) of a rectangle matrix, requiring the same number of colors as the matrix dimensions (respectively 6 and 12 in this case).}
    \label{fig:rectangle}
    \vspace{-15pt}
\end{figure}

\begin{figure}[ht]
    \centering
    \includegraphics[width=0.3\linewidth]{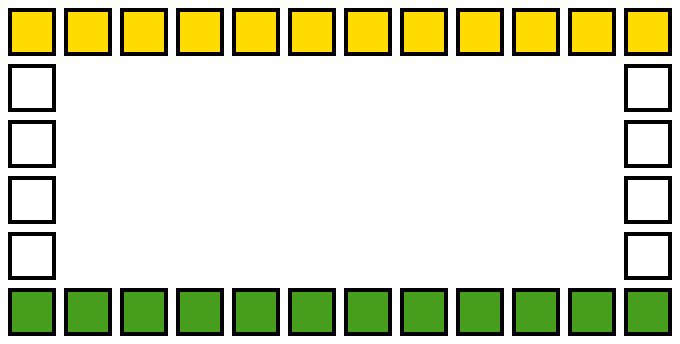}
    \includegraphics[width=0.3\linewidth]{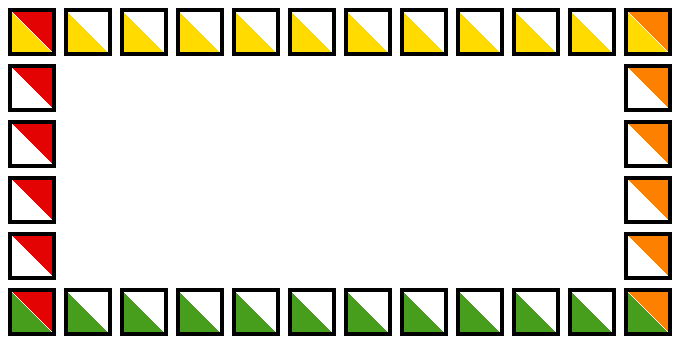}
    \includegraphics[width=0.3\linewidth]{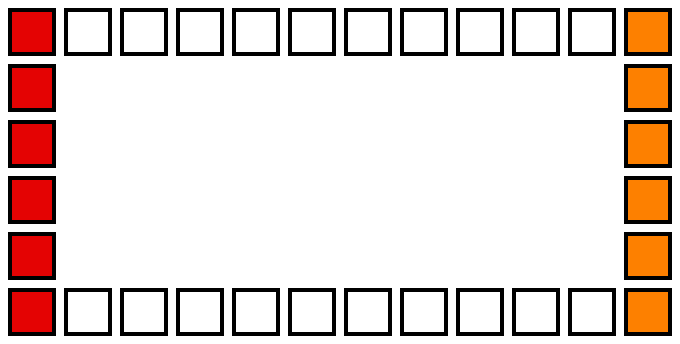}
    \caption{Bicoloring of a rectangle matrix, requiring only 2 colors for the rows (left) and 2 colors for the columns (right). In the central figure, each nonzero coefficient is colored using its row's color and its column's color, when it is not neutral.}
    \label{fig:bicoloring-rectangle}
    \vspace{-15pt}
\end{figure}

In some cases, like for banded matrices, non-overlapping columns and rows are easy to find, and the number of colors remains much lower than the dimension (scaling with the number of bands).
But when a few columns or rows exhibit dense structures, the number of colors required can increase significantly. 
This justifies a shift from unidirectional to bidirectional coloring, or \emph{bicoloring} \citep{hossainComputingSparseJacobian1998,colemanEfficientComputationSparse1998}, which handles columns and rows jointly by combining reverse-mode and forward-mode AD.
We illustrate this idea on a rectangle matrix in \Cref{fig:rectangle}, where a dense row and a dense column impose a significant burden on any one-sided coloring scheme.
As seen in \Cref{fig:bicoloring-rectangle}, bicoloring offers a solution to this problem by assigning colors to a subset of both columns and rows: colored columns to be computed in forward mode and colored rows to be computed in reverse mode.

\subsection{Applications}

In certain optimization problems, the structure of the Jacobian can make unidirectional coloring inefficient.
This subsection highlights cases where bidirectional coloring leads to substantial improvements in the computation of Jacobians.

\subsubsection{Nonlinear least-squares problems}

In the Gauss-Newton method~\citep{bjorckNumericalMethodsLeast2024} for solving $\min_{x \in \mathbb{R}^n} \tfrac{1}{2}\|F(x)\|^2$,
with $F : \mathbb{R}^n \rightarrow \mathbb{R}^m$, one must solve a sequence of linear least-squares subproblems.
At iteration $k$, the subproblem is
$$
    \min_{d \in \mathbb{R}^n}~~\tfrac{1}{2}\|J(x_k)d + F(x_k)\|^2,
$$
where $J(x) \in \mathbb{R}^{m \times n}$ is the Jacobian of $F$ at $x$.
When $m \gg n$, column coloring with forward-mode AD is typically used for efficient Jacobian computation.
However, the presence of a dense row, due to a residual involving all variables (such as normalization constraints or global coupling terms) renders column coloring inefficient, as the entire row must be recovered regardless of compression.
While row coloring with reverse-mode AD is an alternative, it tends to be inefficient when $m$ is much larger than $n$ (due to a large number of row colors).
Bicoloring can improve performance by recovering sparse columns using forward-mode AD and a few dense rows using reverse-mode AD.

\subsubsection{Equality-constrained optimization problems}

For optimization problems of the form
$$
    \min_{x \in \mathbb{R}^n} \phi(x) \quad \text{subject to} \quad c(x) = 0,
$$
where $\phi : \mathbb{R}^n \rightarrow \mathbb{R}$ and $c : \mathbb{R}^n \rightarrow \mathbb{R}^m$, the Jacobian $J(x) \in \mathbb{R}^{m \times n}$ of the constraints $c$ is essential for forming the sequence of linear systems encountered in Newton-based interior-point methods~\citep{wrightPrimaldualInteriorPointMethods1997}.
In particular, the linear system at iteration $k$ is given by the Karush–Kuhn–Tucker conditions:
$$
    \begin{bmatrix}
        \nabla^2_{xx} \mathcal{L}(x_k,y_k) & J(x_k)^{\!\top} \\
        J(x_k)                             & 0
    \end{bmatrix}
    \begin{bmatrix}
        \Delta x \\ \Delta y
    \end{bmatrix}
    =
    -\begin{bmatrix}
        \nabla_x \mathcal{L}(x_k,y_k) \\ c(x_k)
    \end{bmatrix},
$$
where $\mathcal{L}(x,y) = \phi(x) - y^{\top} c(x)$ and $y \in \mathbb{R}^{m}$ is a vector of Lagrange multipliers.
When $m \ll n$, reverse-mode AD is typically more efficient for computing $J(x)$ because it propagates derivatives from a small number of outputs.
Consequently, using a row coloring strategy is appropriate in this setting to maximize the benefits of reverse-mode AD.
However, if a dense column appears in the Jacobian, due to multiple constraints depending on a common variable, a pure row coloring strategy becomes inefficient.
Although column coloring using forward-mode AD is an alternative, bicoloring can improve performance by recovering dense columns using forward-mode AD and sparse rows using reverse-mode AD.

\subsubsection{Optimal control problems}

Discretized optimal control problems~\citep{bettsPracticalMethodsOptimal2020} often yield Jacobians with structured sparsity patterns, including both dense rows and columns.
Consider the following optimal control problem:
$$
    \min_{x(\cdot),\, u(\cdot),\, p} C(x(\cdot), u(\cdot), p)
    \quad \text{subject to} \quad \left\vert ~
    \begin{aligned}
         & \dot{x}(t) = f(x(t), u(t), p), \quad t \in [0, 1], \\
         & \int_{0}^{1} g(x(t), u(t), p)\, dt = 0.
    \end{aligned} \right.
$$
where $C$ is a cost function, $x(t)$ is the state, $u(t)$ is the control, and $p$ represents the parameters that influence the dynamics.
In this formulation, $p$ is also treated as a decision variable and optimized jointly with the state and control trajectories; this is particularly useful when $p$ includes physical parameters or the final time $t_f$ (after an appropriate time rescaling).
Discretizing the interval $t = [0,1]$ using a grid $t_0, t_1, \dots, t_N$ with step sizes $h_i = t_{i+1} - t_i$, we define the discrete variables $X = (x_0, \dots, x_N)$ and $U = (u_0, \dots, u_N)$, where $x_i \approx x(t_i)$ and $u_i \approx u(t_i)$.
Using a Crank–Nicolson (trapezoidal) scheme for the dynamics, we obtain for $i=0,\ldots, N-1$:
$x_{i+1} - x_i - \frac{h_i}{2} \Big( f(x_i, u_i, p) + f(x_{i+1}, u_{i+1}, p) \Big) = 0$.
Since $p$ appears in every dynamic equation, its derivative contributes to every row of the Jacobian, resulting in a dense column.
Similarly, discretizing the integral constraint via the trapezoidal rule gives:
$\sum_{i=0}^{N-1} \frac{h_i}{2} \Big( g(x_i, u_i, p) + g(x_{i+1}, u_{i+1}, p) \Big) = 0$.
This equation aggregates information over the entire time horizon, so its derivative with respect to all decision variables produces a dense row in the Jacobian.
Note that this behavior is independent of the discretization scheme.
Although the Crank–Nicolson scheme is easiest to describe, other schemes exhibit the same behavior.

A common strategy to mitigate these dense structures is to augment the state.
For instance, instead of treating $p$ as an explicit optimization variable, one can introduce an additional state variable $q(t)$ governed by $\dot{q}(t) = 0$, with $q(0) = p$.
Likewise, the integral constraint can be reformulated by introducing an extra state $y(t)$ governed by $\dot{y}(t) = g(x(t), u(t), p)$.
Instead of adding new variables and constraints through these reformulations, bicoloring enables efficient Jacobian computation without increasing the problem size.

\subsection{Related works}

The first approaches for efficient computation of sparse Jacobians and Hessians were proposed by \citet{curtisEstimationSparseJacobian1974} and \citet{powellEstimationSparseHessian1979}, respectively.
Their connection to graph coloring was discovered by
\citet{colemanEstimationSparseJacobian1983},
\citet{colemanEstimationSparseHessian1984}, and
\citet{mccormickOptimalApproximationSparse1983}, 
triggering an ambitious research effort to design coloring algorithms and related software \cite{colemanSoftwareEstimatingSparse1984}.
Because the relevant graph coloring problems are NP-hard, most of these algorithms are heuristic in nature.
For a detailed account of related literature, see the survey by \citet{gebremedhinWhatColorYour2005} and references therein.
There are alternative methods for sparse AD which do not rely on colorings \citep[Chapter 7]{griewankEvaluatingDerivativesPrinciples2008}, but we will not discuss them here.

Given a matrix $A \in \mathbb{R}^{m \times n}$ to compress, one seeks the smallest matrix of seeds $U \in \{0, 1\}^{n \times c}$ (resp. $V \in \{0, 1\}^{m \times c}$) such that the compressed matrix $B = A U$ (resp. $B = V^{\top\!\!} A$) contains all the information needed to recover $A$.
In what follows, this recovery process is also called \emph{decompression}.
For \emph{direct} decompression, the non-zero entries of $A$ should be directly readable in the compressed matrix $B$.
For decompression by \emph{substitution}, the non-zero entries should be retrievable by solving a triangular system of equations.
Choosing the right seeds is equivalent to assigning colors to columns or rows, and then partitioning basis vectors by color.
This coloring interpretation requires encoding the matrix $A$ as a graph: two popular representations are the bipartite graph $\mathcal{G}_b$ linking rows to columns (for general matrices) and the adjacency graph $\mathcal{G}_a$ linking columns with each other (for symmetric matrices).
The requirement for injective decompression define the families of colorings that can be used on these graphs (usually more restrictive than the standard distance-1 coloring).
\Cref{tab:coloringProblemsOverview} summarizes the coloring models used in derivative computation via direct and substitution-based methods. 
In all cases, the lower the number of distinct colors, the more efficient the compression.

\begin{table}[ht]
\centering
\begin{tabular}{l|l|l}
                & \textbf{Direct method} & \textbf{Substitution method}  \\ \hline
 Jacobian, unidirectional & \textit{Distance-2 coloring} & NA \\
 Hessian, unidirectional  &  \textit{Star coloring} & \textit{Acyclic coloring} \\
 Jacobian, bidirectional  & \textit{Star bicoloring} & \textit{Acyclic bicoloring} \\
\end{tabular}
\caption{Overview of coloring problems in derivative computation via direct and substitution methods.}
\label{tab:coloringProblemsOverview}
\vspace{-20pt}
\end{table}

\begin{figure}[htb]
    \centering
    \begin{subfigure}{0.48\textwidth}
    \centering
    \includegraphics[width=\textwidth]{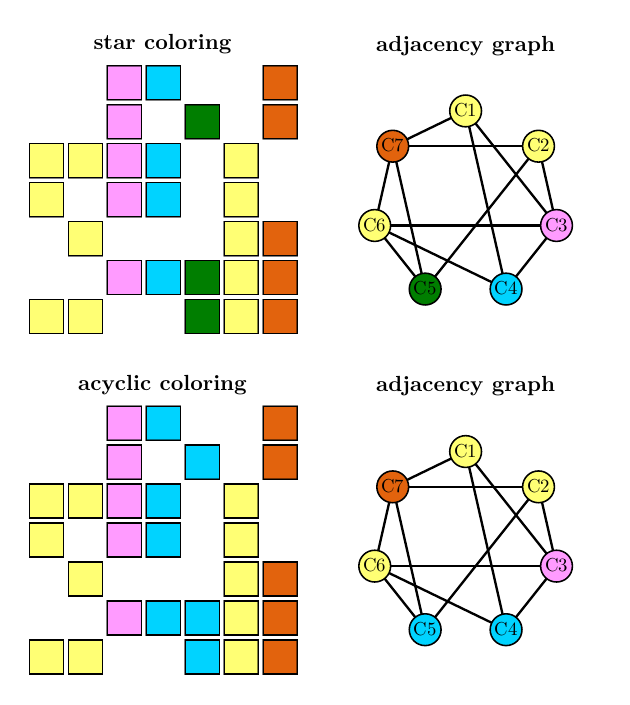}
    \caption{Star and acyclic colorings of a symmetric matrix.}
    \label{fig:graph_symcoloring}
    \end{subfigure}
    \hfill
    \begin{subfigure}{0.48\textwidth}
    \centering
    \includegraphics[width=\textwidth]{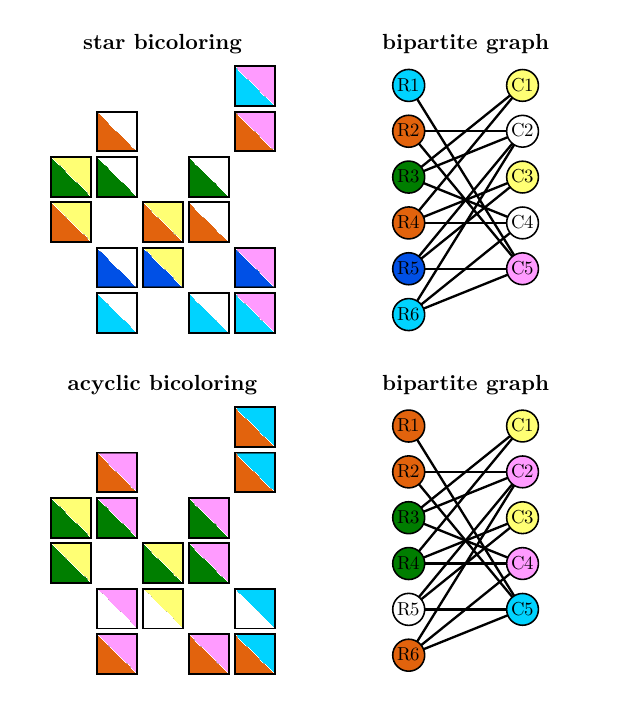}
    \caption{Star and acyclic bicolorings of a non-symmetric matrix.}
    \label{fig:graph_bicoloring}
    \end{subfigure}
    \vspace{-5pt}
    \caption{Example colorings and bicolorings with their graph representations.}
    \vspace{-10pt}
\end{figure}

Central to our work are the studies focused on symmetric matrices (like Hessians).
They give rise to specialized coloring variants: the \emph{star coloring} and \emph{acyclic coloring}, introduced by \citet{colemanEstimationSparseHessian1984} and \citet{colemanCyclicColoringProblem1986} under slightly different names.
We use the terminology of \citet{gebremedhinWhatColorYour2005}, which underlines the importance of \emph{two-colored substructures} in star and acyclic colorings of the adjacency graph $\mathcal{G}_a$.
In a star coloring (which corresponds to direct decompression), given any pair of colors, the induced subgraph is a \emph{collection of stars}.
In an acyclic coloring (which corresponds to decompression by substitution), it is a \emph{collection of trees}.
These insights were exploited by \citet{gebremedhinNewAcyclicStar2007,gebremedhinEfficientComputationSparse2009} to design efficient coloring and decompression routines.
\Cref{fig:graph_symcoloring} showcases examples of star and acyclic colorings.

A coloring model of particular interest to us is bicoloring, which is used in computations of Jacobians and was suggested independently by \citet{hossainComputingSparseJacobian1998} and \citet{colemanEfficientComputationSparse1998}.
The core idea is that combining row and column compression can yield fewer colors overall than either one separately.
In this setting, recovery is based on a pair $(B_c = AU,~B_r = V^{\top\!\!} A)$ of column- and row-compressed matrices rather than a single compressed matrix.
The typical objective is to minimize the total number of seeds in $(U, ~V)$.
Bicoloring may leave some columns (resp. rows) uncolored if all their non-zeros can already be recovered from the corresponding rows (resp. columns).
In such cases, we say that uncolored rows and columns use the \emph{neutral color}.
Examples of bidirectional colorings are displayed on \Cref{fig:graph_bicoloring}.

With direct decompression for bidirectional Jacobian computation in mind, \citet{hossainComputingSparseJacobian1998} elicit the constraint that every path on four vertices in the bipartite graph $\mathcal{G}_b$ should use at least three colors, and deduce a simple greedy algorithm for bicoloring.
\citet{colemanEfficientComputationSparse1998} offer a more general treatment which also applies to decompression by substitution, by specifying that every cycle in $\mathcal{G}_b$ should use at least three colors.
Their algorithm includes a preprocessing step to decide which vertices could be left uncolored, which can be cast as an independent set problem.
\citet{gebremedhinWhatColorYour2005} compare both approaches and suggests connections with the theoretical framework developed for symmetric matrices.
\citet{goyalBiDirectionalDeterminationSparse2006} suggest an Integer Linear Programming (ILP) formulation.
\citet{juedesColoringJacobiansRevisited2012} propose a bicoloring algorithm with a guaranteed approximation ratio on the number of colors.
\citet{hossainGraphModelsTheir2013} explore a new type of graph encoding which they claim is better suited to bidirectional approaches, but does not exhibit invariance by permutation.
\citet{xuEfficientPartialDetermination2013} generalize existing methods to account for other kinds of structure in the Jacobian (like repetition), as well as partial determination.
\citet{gaurDeterminingSparseJacobian2016} give a new lower bound on the number of colors necessary for bicoloring and outline an intricate heuristic algorithm.
\citet{juedesGenericFrameworkApproximation2021} extend their previous study of approximation ratios and introduce new algorithms with performance guarantees.
Finally, a recent paper by \citet{gaurStarBicolouringBipartite2024} leverages a new ILP formulation with column generation to obtain good-quality solutions.

In terms of existing software implementations, \texttt{ColPack} \citep{gebremedhinColPackSoftwareGraph2013} stands out because it was designed as a standalone matrix coloring library in C++, encompassing a wide variety of algorithms.
That is why we choose it as our point of comparison.
To the best of our knowledge, coloring tools are most often part of larger AD frameworks, integrated to enable sparse Jacobians and Hessians.
\texttt{ColPack} itself interfaces with C/C++ AD toolkits like \texttt{ADOL-C} \citep{griewankAlgorithm755ADOLC1996,waltherGettingStartedADOLC2009}, \texttt{CppAD} \citep{bellCppADPackageAlgorithmic2012} or \texttt{ADIC2} \citep{narayananADIC2DevelopmentComponent2010,narayananSparseJacobianComputation2011}.
Sparse matrix coloring implementations for other languages include \texttt{ADMIT} \citep{colemanADMIT1AutomaticDifferentiation2000}, \texttt{ADiMAT} \citep{willkommNewUserInterface2014} and \texttt{MAD} \citep{forthEfficientOverloadedImplementation2006} in MATLAB, \texttt{SparseDiffTools.jl} \citep{juliadiffcontributorsJuliaDiffSparseDiffToolsjl2024,gowdaSparsityProgrammingAutomated2019} in Julia, \texttt{SparseHessianFD} \citep{braunSparseHessianFDPackageEstimating2017} in R, \texttt{sparsejac} \citep{schubertMfschubertSparsejac2024} in Python, as well as the domain-specific languages \texttt{CasADi} \citep{anderssonCasADiSoftwareFramework2019} and \texttt{JuMP.jl} \citep{dunningJuMPModelingLanguage2017,lubinJuMP10Recent2023} for mathematical programming.

\subsection{Contributions}

This paper presents novel bicoloring algorithms, both with direct decompression and decompression via substitution.
In contrast to previous approaches, our work builds on existing algorithms for star and acyclic 
coloring to address the bicoloring problem.
Additionally, we propose a \emph{post-processing} routine for symmetric colorings and bicolorings that aims to neutralize unneeded colors.
In this context, we clarify the relationship between the neutral color and the use of two-colored stars and trees during the decompression of star and acyclic colorings.
Implementation improvements are also described for star and acyclic coloring, as well as for decompression and vertex ordering.

We further introduce the Julia package \texttt{SparseMatrixColorings.jl}\footnote{\texttt{\url{https://github.com/gdalle/SparseMatrixColorings.jl}}}, which provides implementations of our new methods, along with improved versions of classical techniques for row, column, star, and acyclic colorings, as well as ordering and decompression routines. 
It was designed as an alternative to \texttt{ColPack} \citep{gebremedhinColPackSoftwareGraph2013}, written in a high-level programming language \citep{bezansonJuliaFreshApproach2017} with much fewer lines of code but nearly equivalent functionality and comparable speed.

\section{Relation between bidirectional colorings and symmetric colorings} \label{sec:bitosym}

Unlike unidirectional colorings, where an entry $J_{ij}$ can only be directly recovered from the compressed row of its row's color (in row coloring) or the compressed column of its column's color (in column coloring), bidirectional colorings offer more flexibility by leveraging both colors.
This enables more efficient recovery, both directly and through substitution.
This flexibility is also present in symmetric colorings, where any entry $H_{ij}$ in the lower triangular part ($i \ge j$) can be recovered from its symmetric counterpart $H_{ji}$ in the upper triangular part.
The symmetry relation $H_{ij} = H_{ji}$ thus allows leveraging color choices along both the row and column dimensions.

\begin{figure}[ht!]
    \centering
    \begin{subfigure}{0.48\textwidth}
        \centering
        \includegraphics[height=0.45\linewidth]{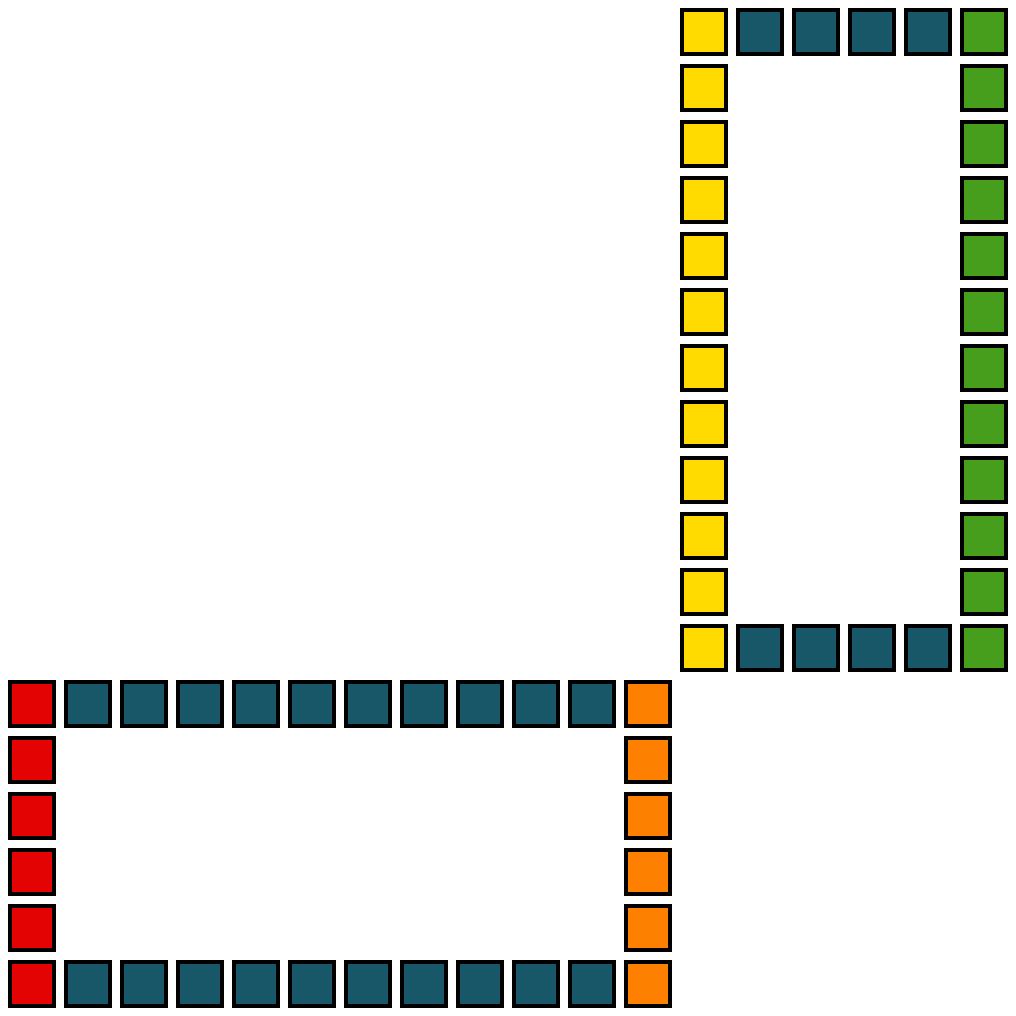}
        \includegraphics[height=0.45\linewidth]{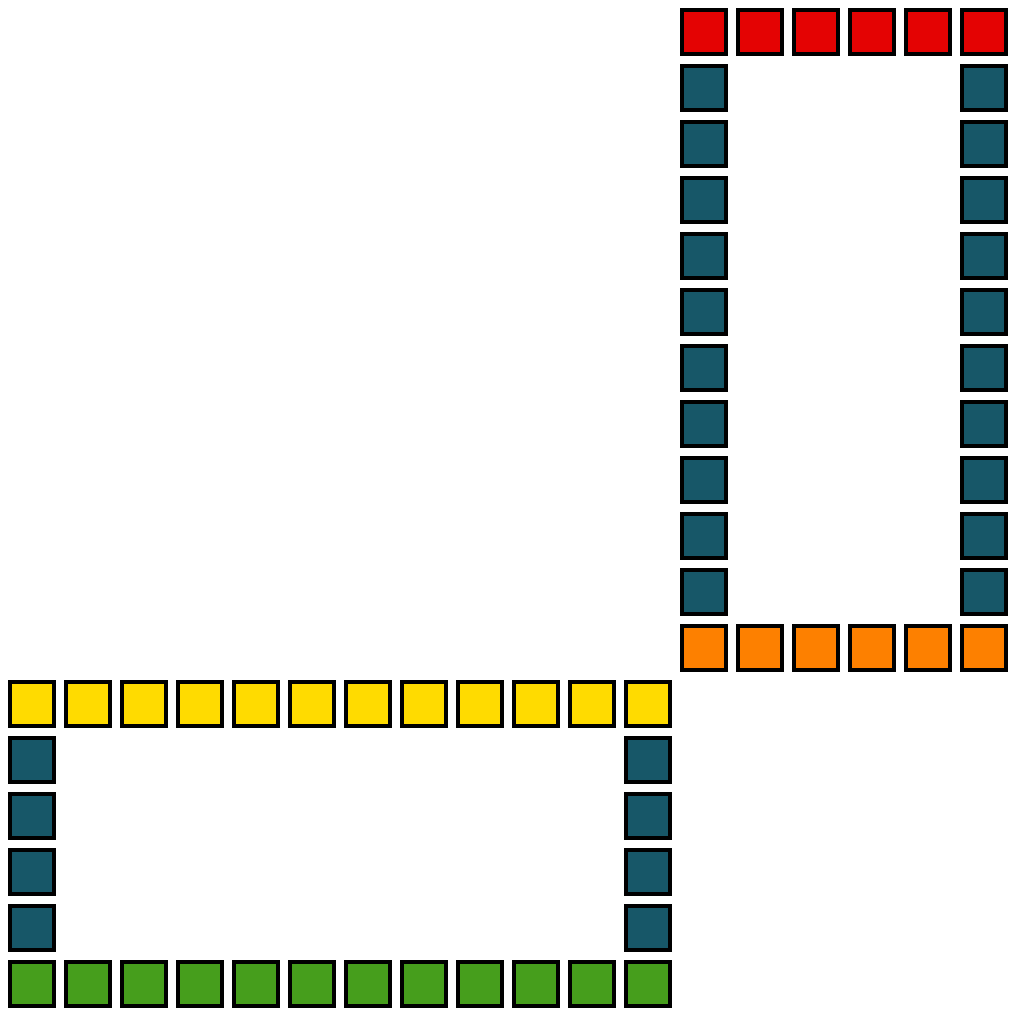}
        \hfill
        \caption{Symmetric coloring on $H$. Nonzeros are colored by the color of their columns on the left panel and by the color of their rows on the right panel.}
        \label{fig:rectangle-H}
        \vspace{-10pt}
    \end{subfigure}
    \hfill
    \begin{subfigure}{0.46\textwidth}
        \centering
        \begin{tabular}{c}
            \includegraphics[width=0.4\linewidth]{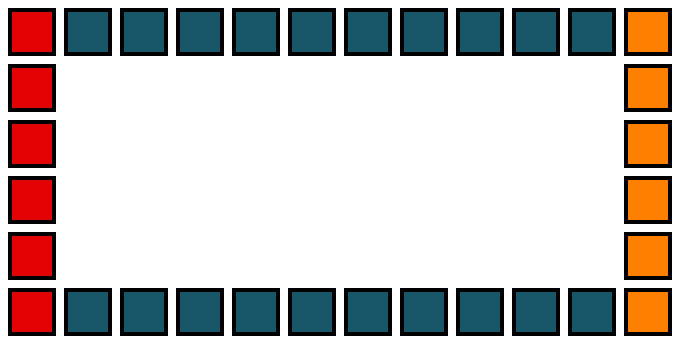} \\
            \includegraphics[width=0.4\linewidth]{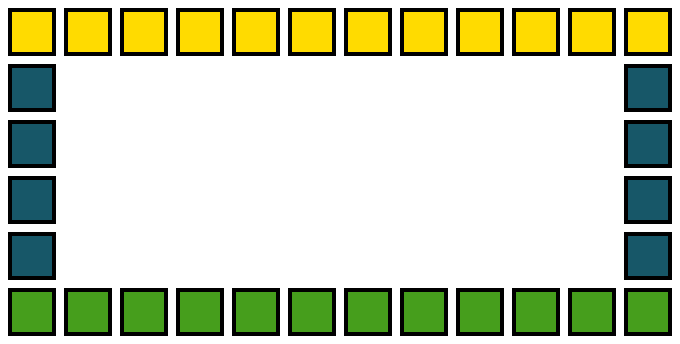}
        \end{tabular}
        \caption{Bicoloring on $J$. Nonzeros are colored according to their column colors in the top panel and according to their row colors in the bottom panel.}
        \label{fig:rectangle-J}
        \vspace{-10pt}
    \end{subfigure}
    \caption{Link between symmetric coloring and bicoloring}
\end{figure}

Building on this observation, we propose a new bidirectional coloring approach that reformulates the bicoloring problem as a symmetric coloring problem.
To achieve this, we construct an \emph{augmented matrix} $H$, where $J$ is stored in the lower triangular portion and $J^\top$ in the upper triangular portion.
This formulation is given by:
\begin{equation}
    H := \begin{bmatrix} ~0 & J^\top~ \\ ~J & 0\phantom{^T}~ \end{bmatrix}.
    \label{eq:augmented_matrix}
\end{equation}
The adjacency graph $\mathcal{G}_a$ of the augmented matrix $H$ is exactly the same as the bipartite graph $\mathcal{G}_b$ of the initial matrix $J$.
With this representation, the symmetric colors assigned to the first $n$ columns (or rows) of $H$ determine the column colors of $J$, while those assigned to the last $m$ columns (or rows) of $H$ determine the row colors of $J$.
Recovering $J$ only requires decompressing the lower triangular portion of $H$, which is a standard procedure for symmetric matrices where only one triangle needs to be stored.
This transformation allows us to extend the well-established symmetric coloring algorithms -- such as star and acyclic coloring techniques -- along with their associated decompression strategies and complexity properties developed by \citet{gebremedhinNewAcyclicStar2007,gebremedhinEfficientComputationSparse2009} to the context of bicoloring.
\Cref{fig:rectangle-H,fig:rectangle-J} illustrate how a symmetric coloring of $H$ translates into a bidirectional coloring of $J$, which exhibits a rectangular sparsity pattern.
Note that a representation of adjoint graph in the context of hypergraphs~\citep{liuNWHyFrameworkHypergraph2022} shares similarities with the structure of the augmented matrix~\eqref{eq:augmented_matrix}.

To preserve the mapping between the symmetric colors of $H$ (stored as an integer vector \texttt{sym\_colors}) and the row and column colors of $J$ (stored as integer vectors \texttt{row\_colors} and \texttt{col\_colors}), we introduce two additional integer vectors: \texttt{sym\_to\_col} and \texttt{sym\_to\_row}.
These vectors associate each symmetric color (numbered from \texttt{1} to \texttt{num\_sym\_colors}) with its corresponding column and row colors.
If an entry in either \texttt{sym\_to\_col} (or \texttt{sym\_to\_row}) is \texttt{0}, it signifies that the corresponding symmetric color is exclusively assigned to a subset of rows (or columns) in $J$.
Additionally, we require that the row and column colors be consecutively numbered starting from \texttt{1}.
This contiguous indexing is crucial because each row or column index in the compressed matrices is directly related to a color.
Maintaining this mapping simplifies the selection of seeds for AD and streamlines the decompression step.
\Cref{alg:remap_colors} illustrates the algorithm used to compute the correspondence between the symmetric colors of $H$ and the row and column colors of $J$.
\begin{algorithm}
    \footnotesize
    \begin{framed}
        \textbf{Input:}
\begin{itemize}
    \setlength{\itemindent}{-23.5pt}
    \item Two integers \texttt{m} and \texttt{n} corresponding to the number of rows and columns of $J$.
    \item An integer \texttt{num\_sym\_colors} representing the total number of symmetric colors.
    \item An integer array \texttt{sym\_colors[1\ldots(n+m)]} specifying a star or acyclic coloring.
\end{itemize}

\textbf{Output:}
\begin{itemize}
    \setlength{\itemindent}{-23.5pt}
    \item \texttt{row\_colors[1\ldots m]}: renumbered colors for rows.
    \item \texttt{col\_colors[1\ldots n]}: renumbered colors for columns.
    \item \texttt{sym\_to\_row[1\ldots num\_sym\_colors]}: maps symmetric colors to row colors.
    \item \texttt{sym\_to\_col[1\ldots num\_sym\_colors]}: maps symmetric colors to column colors.
    \item \texttt{num\_row\_colors}: number of colors for rows.
    \item \texttt{num\_col\_colors}: number of colors for columns.
\end{itemize}

\begin{algorithmic}[1]
    \State \texttt{sym\_to\_col}[1\ldots\texttt{num\_sym\_colors}] $\gets 0$
    \State \texttt{col\_colors[1\ldots n]} $\gets 0$
    \State \texttt{num\_col\_colors} $\gets 0$
    \For{$j \leftarrow 1, \dots, n$}
    \State Let $c_j \leftarrow$ \texttt{sym\_colors[j]}
    \If{$c_j > 0$}
    \If{\texttt{sym\_to\_col[$c_j$]} \texttt{==} 0}
    \State \texttt{num\_col\_colors} $\leftarrow$ \texttt{num\_col\_colors + 1}
    \State \texttt{sym\_to\_col[$c_j$]} $\leftarrow$ \texttt{num\_col\_colors}
    \EndIf
    \State \texttt{col\_colors[j]} $\leftarrow$ \texttt{sym\_to\_col[$c_j$]}
    \EndIf
    \EndFor
    \State \texttt{sym\_to\_row[1\ldots num\_sym\_colors]} $\gets 0$
    \State \texttt{row\_colors[1\ldots m]} $\gets 0$
    \State \texttt{num\_row\_colors} $\gets 0$
    \For{$i \leftarrow$ \texttt{n+1,\dots,n+m}}
    \State Let $c_i \leftarrow$ \texttt{sym\_colors[i]}
    \If{$c_i > 0$}
    \If{\texttt{sym\_to\_row[$c_i$]} \texttt{==} 0}
    \State \texttt{num\_row\_colors} $\leftarrow$ \texttt{num\_row\_colors + 1}
    \State \texttt{sym\_to\_row[$c_i$]} $\leftarrow$ \texttt{num\_row\_colors}
    \EndIf
    \State \texttt{row\_colors[i-n]} $\leftarrow$ \texttt{sym\_to\_row[$c_i$]}
    \EndIf
    \EndFor
\end{algorithmic}

    \end{framed}
    \caption{Procedure \texttt{remap\_colors} for mapping symmetric colors to row and column colors. The arrays use 1-based indexing.}
    \label{alg:remap_colors}
    \vspace{4pt}
\end{algorithm}

Since our approach relies on symmetric colorings, the decompression routines require access to the compressed symmetric matrix $B_s$ (or its transpose).
$B_s$ corresponds to the compression of the columns of $H$, grouped according to symmetric colors.
Compressing by rows instead yields $B_s^\top$.
Thanks to our mapping, explicitly constructing $B_s$ is unnecessary.
Each coefficient of $B_s$ can be efficiently retrieved from $B_r$ and $B_c$, where $B_r$ stores the compressed rows of $J$, grouped by row colors, and $B_c$ stores the compressed columns of $J$, grouped by column colors.
The following holds:
\begin{equation}
\resizebox{0.9\textwidth}{!}{$
    B_s\texttt{[r,c]} =
    \begin{cases}
        B_r\texttt{[\texttt{sym\_to\_row[c]}, r]} \!\!   & \text{if } 1 \leq \texttt{r} \leq n \text{ and } \texttt{sym\_to\_row[c]} \neq 0,   \\[1mm]
        B_c\texttt{[r-n, \texttt{sym\_to\_col[c]}]} \!\! & \text{if } 1 \leq \texttt{r}-n \leq m \text{ and } \texttt{sym\_to\_col[c]} \neq 0, \\[1mm]
        \texttt{0} \!\!                                  & \text{otherwise.}
    \end{cases}
$}
\label{eq:Bs_Br_Bc}
\end{equation}
Bands of zeros in $B_s$ that are not needed during decompression are not explicitly stored in $B_r$ and $B_c$.
They correspond to the structural zero blocks in $H$.

In \Cref{eq:augmented_matrix}, we could have chosen to put $J^\top$ in the lower triangular portion and $J$ in the upper triangular portion instead.
This switch is equally valid, and only requires minor adaptations in our algorithms.
It would amount to reverting the order of vertices in the augmented adjacency graph $\mathcal{G}_a$ obtained from the bipartite graph $\mathcal{G}_b$, listing row vertices first and column vertices next.

\section{Post-processing} \label{sec:postprocessing}

\subsection{Neutral color in symmetric and bidirectional colorings}

We can see in \Cref{fig:rectangle-H,fig:rectangle-J} that the dark blue color is unnecessary, so the rows and columns currently using it can be left uncolored instead, as illustrated in \Cref{fig:bicoloring-rectangle}.
A neutral color indicates that the coefficients of neutralized columns or rows are not needed to recover all non-zeros.
This is materialized by a blank square or triangle in our visualizations.

We developed a post‐processing algorithm to reduce the number of colors by replacing some of them with a neutral color after applying symmetric colorings, thereby also reducing the number of colors in our bicolorings.
Thanks to symmetry, each off-diagonal entry can be recovered from either its row or its column, which are associated with different indices.
In previous symmetric coloring approaches, every diagonal coefficient was assumed to be nonzero, making each color essential for diagonal retrieval and leaving no room for neutral colors.
In contrast, for our bicoloring approach, all diagonal elements of $H$ are zero, making post-processing beneficial.
Note that our approach is not specific to the $2 \times 2$ block structure of $H$ used for bicoloring.
It can be applied to any symmetric matrix with at least some zeros on its diagonal and that has been colored using star or acyclic coloring.

After our post-processing, the non‑neutral (i.e.\ colored) row and column vertices constitute a \emph{vertex cover} of the corresponding graph ($\mathcal{G}_a$ or $\mathcal{G}_b$).
Indeed, every nonzero entry (edge in the graph) is incident to at least one colored vertex.

\subsection{Relation between neutral color and two-colored structures}

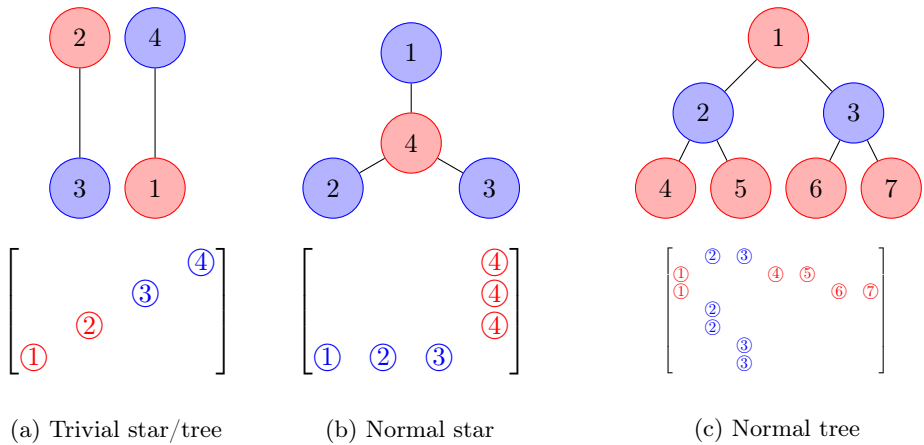
\begin{figure}[ht!]
    \begin{subfigure}[t]{0.25\textwidth}
        \centering
        \begin{tikzpicture}
            \tikzstyle{red_node} = [circle, draw=red, fill=red!30, minimum size=8mm, inner sep=0pt]
            \tikzstyle{red_node2} = [circle, draw=red, fill=red!30, minimum size=8mm, inner sep=0pt]
            \tikzstyle{blue_node} = [circle, draw=blue, fill=blue!30, minimum size=8mm, inner sep=0pt]
            \tikzstyle{blue_node2} = [circle, draw=blue, fill=blue!30, minimum size=8mm, inner sep=0pt]

            \node[red_node] (t1) at (1,0) {1};
            \node[red_node2] (t2) at (0,2) {2};
            \node[blue_node] (t4) at (1,2) {4};
            \node[blue_node2] (t3) at (0,0) {3};
            \draw (t1) -- (t4);
            \draw (t2) -- (t3);
        \end{tikzpicture}
        $$
        \bmat{
           & & & \textcolor{blue}{\circled{4}}
        \\ & & \textcolor{blue}{\circled{3}} &
        \\ & \textcolor{red}{\circled{2}} & &
        \\ \textcolor{red}{\circled{1}} & & &
        }
        $$
        \caption{Trivial star/tree}
        \label{fig:2trivialstar}
    \end{subfigure}
    \hfill
    \begin{subfigure}[t]{0.3\textwidth}
        \centering
        \begin{tikzpicture}
            \tikzstyle{red_node} = [circle, draw=red, fill=red!30, minimum size=8mm, inner sep=0pt]
            \tikzstyle{blue_node} = [circle, draw=blue, fill=blue!30, minimum size=8mm, inner sep=0pt]

            \node[red_node] (center) at (0,0) {4};

            \foreach \i [count=\j from 1] in {90, 210, 330} {
                    \node[blue_node] (leaf\j) at (\i:1.2) {\j};
                    \draw (center) -- (leaf\j);
                }
        \end{tikzpicture}
        $$
        \bmat{
           &     &     & \textcolor{red}{\circled{4} }
        \\ &     &     & \textcolor{red}{\circled{4} }
        \\ &     &     & \textcolor{red}{\circled{4} }
        \\ \textcolor{blue}{\circled{1}} & \textcolor{blue}{\circled{2}} & \textcolor{blue}{\circled{3}} &
        }
        $$
        \caption{Normal star}
        \label{fig:2normalstar}
    \end{subfigure}
    \hfill
    \begin{subfigure}[t]{0.4\textwidth}
        \centering
        \begin{tikzpicture}
            \tikzstyle{red_node} = [circle, draw=red, fill=red!30, minimum size=8mm, inner sep=0pt]
            \tikzstyle{blue_node} = [circle, draw=blue, fill=blue!30, minimum size=8mm, inner sep=0pt]

            \node[red_node] (root) at (0,0) {1};

            \node[blue_node] (left) at (-1,-1) {2};
            \node[blue_node] (right) at (1,-1) {3};
            \draw (root) -- (left);
            \draw (root) -- (right);

            \node[red_node] (leftleft) at (-1.5,-2) {4};
            \node[red_node] (leftright) at (-0.5,-2) {5};
            \node[red_node] (rightleft) at (0.5,-2) {6};
            \node[red_node] (rightright) at (1.5,-2) {7};

            \draw (left) -- (leftleft);
            \draw (left) -- (leftright);
            \draw (right) -- (rightleft);
            \draw (right) -- (rightright);
        \end{tikzpicture}
        $$
        \resizebox{0.59\textwidth}{!}{$
        \bmat{
           & \textcolor{blue}{\circled{2}} & \textcolor{blue}{\circled{3}} & & & &
        \\ \textcolor{red}{\circled{1}} & & & \textcolor{red}{\circled{4}} & \textcolor{red}{\circled{5}} & &
        \\ \textcolor{red}{\circled{1}} & & & & & \textcolor{red}{\circled{6}} & \textcolor{red}{\circled{7}}
        \\ & \textcolor{blue}{\circled{2}} & & & & &
        \\ & \textcolor{blue}{\circled{2}} & & & & &
        \\ & & \textcolor{blue}{\circled{3}} & & & &
        \\ & & \textcolor{blue}{\circled{3}} & & & &
        }
        $
        }
        $$
        \caption{Normal tree}
        \label{fig:2normaltree}
    \end{subfigure}
    \caption{Variants of two-colored structures with example matrices.}
    \label{fig:two-colored-structures}
    \vspace{-10pt}
\end{figure}

In graph theory, star and acyclic colorings are specialized forms of vertex colorings that impose additional constraints to prevent certain substructures within the graph.
A star coloring is a proper vertex coloring where no path on four vertices uses two colors, meaning any path on four vertices uses at least three distinct colors.
This restriction ensures that the subgraph induced by any two color structures forms a disjoint union of star graphs, thereby eliminating bicolored four-vertex paths.
An acyclic coloring is a proper vertex coloring such that every cycle in the graph is colored with at least three colors, ensuring that no two-colored cycles exist. In other words, the subgraph induced by any two color structures is a forest, thereby eliminating bicolored cycles.
\Cref{fig:two-colored-structures} illustrates the variants of two-colored structures encountered in star and acylic colorings.

To determine if a color can be discarded, we need to identify cases where the color is not required during the decompression phase in star and acyclic colorings, which rely on two-colored structures.
The first prerequisite for neutralizing a color in symmetric colorings is that it is not needed to recover diagonal coefficients.
In the context of star coloring, only the color of the hubs of the stars is needed during decompression.
The color of the spokes is not used (see the routine
\textit{DirectRecover2} in~\citet{gebremedhinEfficientComputationSparse2009}).
To illustrate this, consider an adjacency matrix $J$ with one row and three columns.
The associated augmented matrix $H$ is represented on \Cref{fig:2normalstar}):
Here, the columns colored blue (spokes) are not used; only the column colored red (hub) is needed.
Consequently, any color that does not appear in the hubs of any two-colored stars can safely be replaced by the neutral color.

In the context of acyclic coloring, only the color of the leaves in trees is not needed during decompression (see the routine \textit{IndirectRecover} in~\citet{gebremedhinEfficientComputationSparse2009}).
However, for two-colored trees that are neither trivial nor standard two-colored stars (trees with a depth of at least 2), both colors must be retained as they correspond to the root and internal nodes (see the graph in \Cref{fig:2normaltree}).
Therefore, the colors that can be neutralized in trees are those found in edges and stars (trees with a depth of 1).

The main challenge comes from edges, which are simplified forms of two-colored stars and trees that can be found in both star and acyclic colorings.
In these cases, the color of either vertex can be used for decompression.
For that reason, we suggest first processing all two-colored stars and trees with at least three vertices.
Then, at the end of the procedure, we iterate over the remaining edges, leveraging our knowledge of nontrivial stars and trees.
For each edge, if either of its two vertices is already (a) a hub of a normal two-colored star or (b) is simply a vertex in a normal two-colored tree, we cannot give it the neutral color.
Thus, we choose to use that vertex for decompression of the current edge.
Therein lies the benefit of processing individual edges at the end: if one color is already indispensable, we pick it for the hub, instead of making the other color indispensable too when it is not needed.
On the other hand, if both colors can still be neutralized in the edge, we must arbitrarily neutralize one and then ensure that decompression uses the color of the other vertex.
It may require updating who the hub is in these trivial two-colored structures.

For example, consider the following anti-diagonal matrix of \Cref{fig:2trivialstar}.
We can remove either the blue or red color to compute all non-zeros.
\Cref{alg:postprocessing} summarizes the post-processing algorithm for star and acyclic colorings.
Note that it can be simplified if we only use it for bicoloring, because $H$ always has a diagonal of zeros, as well as in star coloring, because the set of normal two-colored trees $\mathcal{T}$ is always empty.

\begin{algorithm}[htp!]
    \footnotesize
    \begin{framed}
        \textbf{Input:}
\begin{itemize}
    \setlength{\itemindent}{-23.5pt}
    \item An integer \texttt{p} corresponding to the dimension of a symmetric matrix $H$.
    \item An integer \texttt{num\_sym\_colors} representing the total number of symmetric colors.
    \item An integer array \texttt{sym\_colors[1\ldots p]} specifying a star or acyclic coloring.
    \item An integer array \texttt{offsets[1\ldots p]} used as a workspace.
    \item A boolean vector \texttt{used[1\ldots num\_sym\_colors]}: colors that can be neutralized.
    \item A set $\mathcal{T}$ of normal two-colored trees.
    \item A set $\mathcal{S}$ of normal two-colored stars.
    \item A set $\mathcal{E}$ of trivial two-colored structures.
\end{itemize}

\textbf{Output:}
\begin{itemize}
    \setlength{\itemindent}{-23.5pt}
    \item The integer vector \texttt{sym\_colors[1\ldots p]} updated.
\end{itemize}

\begin{algorithmic}[1]
    \State \texttt{used[1\ldots num\_sym\_colors]} $\gets$ \texttt{false}
    \For{$k \leftarrow 1, \dots, p$}
    \If{$H_{k,k} \ne 0$}
    \State \texttt{used[sym\_colors[k]] = true}
    \EndIf
    \EndFor
    \For{\textbf{each} normal two-colored tree T in $\mathcal{T}$}
    \State Let $i$ and $j$ be the indices of two vertices connected by an edge in T
    \State \texttt{used[sym\_colors[i]] = true}
    \State \texttt{used[sym\_colors[j]] = true}
    \EndFor
    \For{\textbf{each} normal two-colored star S in $\mathcal{S}$}
    \State Let $h$ be the index of the hub vertex in S
    \State \texttt{used[sym\_colors[h]] = true}
    \EndFor
    \For{\textbf{each} trivial two-colored structure E in $\mathcal{E}$}
    \State Let $i$ and $j$ be the indices of the two vertices of the edge E

    \If{\texttt{used[sym\_colors[i]]} \textbf{or} \texttt{used[sym\_colors[j]]}}
    \State Choose the hub $h$ as $i$ or $j$ such that its color is used
    \Else
    \State Choose the hub $h$ arbitrary between $i$ and $j$
    \State \texttt{used[sym\_colors[h]] = true}
    \EndIf
    \EndFor
    \State \texttt{num\_neutral\_colors} $\gets 0$
    \For{\texttt{c} $\gets$ \texttt{1} \textbf{to} \texttt{num\_sym\_colors}}
    \If{(\textbf{not} \texttt{color\_used[c]})}
    \State \texttt{num\_neutral\_colors} $\gets$ \texttt{num\_neutral\_colors + 1}
    \Else
    \State \texttt{offsets[c]} $\gets$ \texttt{num\_neutral\_colors}
    \EndIf
    \EndFor
    \For{\texttt{i} $\gets$ \texttt{1} \textbf{to} \texttt{p}}
    \State \texttt{c} $\gets$ \texttt{sym\_colors[i]}
    \If{(\textbf{not} \texttt{used[c]})}
    \State \texttt{sym\_color[i]} $\gets 0$
    \Else
    \State \texttt{sym\_color[i]} $\gets$ \texttt{sym\_color[i]} $-$ \texttt{offsets[c]}
    \EndIf
    \EndFor
\end{algorithmic}

    \end{framed}
    \caption{Procedure \texttt{post\_processing} for discarding colors after star and acyclic colorings. The arrays use 1-based indexing.}
    \label{alg:postprocessing}
    \vspace{4pt}
\end{algorithm}

Because an arbitrary choice may be necessary to neutralize one color in trivial two-colored structures, the total number of neutralized colors may vary.
In the context of bicoloring, this choice can be made to reduce the number of colors in either the column partition or the row partition of a Jacobian.
Depending on the shape of the Jacobian, as well as the cost of forward and reverse AD, one option may be more advantageous than the other.

\subsection{Post-processing without color elimination}

We recall that the post-processing aims to minimize the total number of colors, rather than maximizing the assignments of the neutral color.
In some cases, the color assignments of rows and columns could be replaced with the neutral color without reducing the total number of distinct colors.
This change affects only the seeds and does not benefit most AD tools.
Tools like \texttt{Tapenade}~\citep{hascoetTapenadeAutomaticDifferentiation2013} utilize activity analysis to identify active variables and skip computations on inactive ones.
While this approach may optimize the process, further research is needed to determine if practical performance gains can be achieved.
Additionally, this optimization could require different code being used to compute $Ju$ and $J^{\top\!} v$, depending on the seeds.

\section{A new kind of coloring}

In the symmetric setting, replacing some colors with the neutral color $0$
means that the end result is no longer a standard star or acyclic coloring.
We shed light on the resulting object in the case of star coloring only, leaving acyclic coloring for future work.

\subsection{Partition perspective}

First, we recall the goal of coloring a symmetric matrix $A$ with direct decompression: obtaining a compressed representation from which every non-zero coefficient can be read directly.
Following \citet[Definition 4.2]{gebremedhinWhatColorYour2005}, a \emph{symmetrically orthogonal partition} of the columns of $A$ is such that for every $A_{ij} \neq 0$, either (1) $A_j$ is the only column in its group with a non-zero entry in row $i$, or (2) $A_i$ is the only column in its group with a non-zero entry in row $j$. 
\citet{colemanEstimationSparseHessian1984} prove that every star coloring of the adjacency graph of $A$ induces a symmetrically orthogonal partition of the columns of $A$.

On the other hand, our post-processing method modifies the star coloring to induce a \emph{neutralized symmetrically orthogonal partition}: it is such that for every $A_{ij} \neq 0$, either (1) $A_j$ has a non-neutral color and is the only column in its group with a non-zero entry in row $i$, or (2) $A_i$ has a non-neutral color and is the only column in its group with a non-zero entry in row $j$.
We now give a graph characterization of this object in the special case where all diagonal coefficients of $A$ are zero, which applies to the augmented matrix of \eqref{eq:augmented_matrix}.
The usual star coloring corresponds to the opposite case where all diagonal coefficients are non-zero, which means the neutral color cannot be used.
Intermediate cases where only some coefficients are non-zero may also be described, but the ensuing definitions are less crisp because the adjacency graph is not sufficient to distinguish between zero and non-zero diagonal entries (unless self-loops are considered).

\subsection{Graph perspective}

Let $A$ be a symmetric $n \times n$ matrix whose diagonal entries are all zero.
Let $\phi : \{1, \dots, n\} \to \{1, \dots, n\}$
be a mapping from the columns of $A$ to integer colors.
Any such mapping defines a partition $P_\phi$ of the columns of $A$, by gathering columns with the same color $c$ inside a single group $P_\phi[c] := \{k \in \{1, \dots, n\}: \phi(k) = c\}$.
As before, $0$ denotes a neutral color, which is not used for recovery of the matrix.
Let $\mathcal{G}_a = (\mathcal{V}, \mathcal{E})$ denote the adjacency graph of $A$, with vertex set $\mathcal{V} = \{1, \dots, n\}$ and edge set $\mathcal{E} = \{(i, j): A_{ij} \neq 0\}$.

\paragraph{Definition}
We say that $\phi$ is a \emph{no-zig-zag coloring} of the adjacency graph $\mathcal{G}_a$ if the following conditions are all met:
\begin{enumerate}
\item Every edge is incident on at least one colored vertex: if $(i, j) \in \mathcal{E}$, then $\phi(i) > 0$ or $\phi(j) > 0$.
\item Neutral-colored vertices have distinctly-colored neighbors: if $j$ and $q$ are adjacent to the same $i$ such that $\phi(i) = 0$, then $\phi(j) \neq \phi(q)$.
\item There is no color-alternated (zig-zagging) path on four vertices: for any path $(q, i, j, p)$ in $\mathcal{G}_a$, either $\phi(q) \neq \phi(j)$ or $\phi(i) \neq \phi(p)$, so that we cannot have a succession of colors of the form $(c_1, c_2, c_1, c_2)$.
\end{enumerate}

Note that this definition does not require $\phi$ to be a coloring of $\mathcal{G}_a$ in the usual sense: adjacent vertices are allowed to share the same non-neutral color $c > 0$.
Therefore, no-zig-zag coloring is a relaxation of star coloring, as introduced by \citet{colemanEstimationSparseHessian1984}.
It is also related to star bicoloring, as introduced by \citet{colemanEfficientComputationSparse1998}.

\paragraph{Theorem}
The mapping $\phi$ is a no-zig-zag coloring of $\mathcal{G}_a$ if and only if the partition $P_\phi$ is a neutralized symmetrically orthogonal partition of the columns of $A$.

Before the formal proof, we give an intuitive way to obtain it.
Given the matrix $A$ below, we check which two-colored partitions allow recovery of $A_{ij}$ (up to color switch).
The results show that the only forbidden coloring pattern for path $(q, i, j, p)$ is the zig-zag.

\begin{minipage}{0.35\textwidth}
    \begin{equation*}
        A = \begin{pmatrix}
        \cdot & A_{qi} & \cdot & \cdot \\
        A_{iq} & \cdot & A_{ij} & \cdot \\
        \cdot & A_{ji} & \cdot & A_{jp} \\
        \cdot & \cdot & A_{pj} & \cdot
        \end{pmatrix}
    \end{equation*}
\end{minipage} \hfill
\begin{minipage}{0.55\textwidth}
\begin{tabular}{lllll}
    \hline
    $q$ & $i$ & $j$ & $p$ & recoverable $A_{ij}$ \\ \hline
    $c_1$ & $c_2$ & $c_1$ & $c_2$ & no \\
    $c_1$ & $c_1$ & $c_2$ & $c_2$ & yes (in $c_1$ or $c_2$) \\
    $c_1$ & $c_2$ & $c_2$ & $c_1$ & yes (in $c_2$) \\
    $c_1$ & $c_1$ & $c_1$ & $c_2$ & yes (in $c_1$) \\
    $c_1$ & $c_1$ & $c_2$ & $c_1$ & yes (in $c_2$) \\ \hline
    \end{tabular}
\end{minipage}

\bigskip

\paragraph{Proof} We tackle the two implications sequentially.
\begin{description}
\item[$\implies$] Let $\phi$ define a no-zig-zag coloring of $\mathcal{G}_a$.
If the partition $P_\phi$ is not a neutralized symmetrically orthogonal for $A$, there exists $(i, j)$ such that $A_{ij} \neq 0$ cannot be recovered from $i$ nor from $j$. This means we are in one of the following situations (the conditions refer to those defining the no-zig-zag coloring):
\begin{itemize}
  \item $\phi(i) = 0$. By condition 1, $\phi(j) > 0$. Since $A_{ij}$ cannot be recovered from $j$, there exists a column $q$ with $A_{iq} \neq 0$ such that $\phi(q) = \phi(j)$. This violates condition 2.
  \item $\phi(j) = 0$. By condition 1, $\phi(i) > 0$. Since $A_{ij}$ cannot be recovered from $i$, there exists a column $p$ with $A_{jp} \neq 0$ such that $\phi(p) = \phi(i)$. This violates condition 2.
  \item $\phi(i) > 0$ and $\phi(j) > 0$. Since $A_{ij}$ cannot be recovered from $i$ nor from $j$, there exist $q$ and $p$ such that $A_{iq} \neq 0$, $A_{jp} \neq 0$, $\phi(i) = \phi(p)$ and $\phi(j) = \phi(q)$. But then the path $(q, i, j, p)$ is color-alternated, which violates condition 3.
\end{itemize}
\item[$\impliedby$] Let $P_\phi$ be a neutralized symmetrically orthogonal partition for $A$. We now check the conditions of the no-zig-zag coloring one by one:
\begin{enumerate}
  \item For $A_{ij} \neq 0$ to be recoverable, either $\phi(i) > 0$ or $\phi(j) > 0$ must hold. So the edge $(i, j) \in \mathcal{E}$ is incident on at least one colored vertex.
  \item If $\phi(i) = 0$ and $i$ has two neighbors $j$ and $q$, then $A_{ij}$ can only be recovered from $j$ and $A_{iq}$ from $q$. This implies that $\phi(j) \neq \phi(q)$.
  \item Let $(q, i, j, p)$ be a path in $\mathcal{G}_a$.
  \begin{itemize}
      \item Assume $\phi(i) = 0$ (resp. $\phi(j) = 0$). Then, condition 2 imposes distinctly colored neighbors for $i$ (resp. $j$). This brings the total number of colors to at least 3 and precludes a color-alternated path.
      \item Now assume $\phi(i) > 0$ and $\phi(j) > 0$. Since $A_{ij}$ is recoverable, this means that either $\phi(p) \neq \phi(i)$ or $\phi(q) \neq \phi(j)$. Thus, the path cannot be color-alternated either.
  \end{itemize}
\end{enumerate}
\end{description}

\emph{In the rest of the paper, even when our star and acyclic coloring algorithms are enhanced with post-processing, we still refer to them under the same name.}
This makes the comparison with previous versions more concise, in particular for benchmark tables.
The reader should keep in mind that the actual object we obtain is a relaxation of the usual symmetric coloring concepts.

\section{Software implementation}

\texttt{SparseMatrixColorings.jl} (or \texttt{SMC} for short) is an open source library containing our new coloring algorithms, as well as some standard ones.
It can be seen as a modernization of \texttt{ColPack}~\citep{gebremedhinColPackSoftwareGraph2013}, written in a high-level language that facilitates applications.
The implementation is based partly on a series of papers by the \texttt{ColPack} authors \citep{gebremedhinWhatColorYour2005,gebremedhinNewAcyclicStar2007,gebremedhinEfficientComputationSparse2009,gebremedhinColPackSoftwareGraph2013}, but includes new algorithmic ideas.
In particular, for bidirectional problems, \texttt{ColPack} only supports star bicoloring, while \texttt{SMC} includes our new acyclic bicoloring.
As the benchmarks will demonstrate, this can lead to a smaller number of colors.

Furthermore, we fixed several bugs in \texttt{ColPack}, developed a C interface, and created a Julia wrapper \texttt{ColPack.jl}\footnote{\texttt{\url{https://github.com/Exanauts/ColPack.jl}}}. The resulting fork\footnote{\texttt{\url{https://github.com/amontoison/ColPack}}} of \texttt{ColPack} 
is made publicly available on GitHub.

\subsection{Features and design}

Like \texttt{ColPack}, our \texttt{SMC} package offers a variety of coloring and decompression options.
It can compute either column, row or bidirectional colorings, for symmetric or non-symmetric matrices, with direct or indirect decompression.
This translates into five different core algorithms: partial distance-2 coloring, star and acyclic coloring, and star and acyclic bicoloring.
All of these algorithms employ greedy procedures that color vertices one after another.
Since the order in which vertices are considered impacts the number of colors used, we allow users to select an ordering criterion:  natural, random, largest first, dynamic largest first, smallest last, or incidence degree \citep{gebremedhinColPackSoftwareGraph2013}.
For each variant of the coloring problem, we only offer the most suitable model and leave alternative approaches aside.
As an example, we provide a partial distance-2 coloring for bipartite graphs, but not distance-1 coloring on intersection graphs (which would be computationally more expensive \citep{gebremedhinWhatColorYour2005}).

Thanks to our Julia implementation, \texttt{SMC} enjoys the benefits of just-in-time compilation and multiple dispatch \citep{bezansonJuliaFreshApproach2017}.
Fast code is generated for a variety of input formats and number types, without needing to anticipate each one.
At the same time, the use of a high-level language enables remarkable conciseness: according to \texttt{cloc}\footnote{\texttt{\url{https://github.com/AlDanial/cloc}}}, the entire package requires around 5 000 lines of Julia code, against nearly 30 000 lines of C\texttt{++} for \texttt{ColPack}, with comparable functionalities.

\subsection{Fast and flexible decompression}
We took great care to optimize the performance of \texttt{SMC}, with a particular focus on decompression.
Indeed, while coloring is usually performed just once, and then amortized over every subsequent Jacobian or Hessian computation, decompression is an integral part of derivative computations and often happens inside hot loops.
Thus, the coloring step in \texttt{SMC} includes pre-allocation of necessary memory and pre-computation of relevant information, to make decompression as efficient as possible.

While \texttt{SMC} allows decompression into arbitrary matrices, it is most efficient when presented with sparse storage.
Typical formats for sparse matrices include Compressed Sparse Column (CSC) or Compressed Sparse Row (CSR).
The CSC format, which Julia considers the default, comprises three arrays:
\begin{itemize}[leftmargin=0pt,labelsep=5pt]
    \setlength{\itemindent}{10pt}
    \item \texttt{nzval}, containing the nonzero values of the matrix, concatenated column‑wise;
    \item \texttt{rowval}, containing the corresponding row indices, sorted in ascending order;
    \item \texttt{colptr}, an integer array of length $n+1$ such that, for each column~\texttt{j}, its nonzero entries appear in \texttt{nzval} and \texttt{rowval} at positions \texttt{colptr[j]} through \texttt{colptr[j+1]-1}.
\end{itemize}
Accessing a single coefficient \texttt{A[i,j]} in a CSC matrix \texttt{A} is expensive.
 It requires a binary search to find \texttt{k} on the sorted slice \texttt{colptr[j]:colptr[j+1]-1} such that \texttt{rowval[k] = i}, after which the requested value is read from or written to at \texttt{nzval[k]}.
Thus, a naive unidirectional decompression \texttt{A[i,j] = B[i,color[j]]} will be slower than it needs to be.
That is why we precompute a vector of compressed indices such that \texttt{nzval[k] = B[compressed\_indices[k]]}, where the matrix \texttt{B} is indexed linearly columnwise.
With this preliminary step out of the way, decompression turns into a simple loop over all \texttt{k in 1:nnz(A)} without binary search.
The computation of \texttt{compressed\_indices} is a one-time cost and is quickly amortized when we need to perform multiple decompressions, as is the case in non-linear optimization methods that require many Jacobians and Hessians.
It also allows decompression on SIMD architectures such as GPUs.

Anticipating cases where $B$ itself might be too expensive to store, we also provide utilities for decompressing a single color at a time when it is feasible (which amounts to $B$ having a single column or row).
Similarly, in the case of symmetric matrices, we offer decompression inside a single triangle (lower or upper), thus skipping duplicate work.
As noted earlier, decompression into a triangle for star and acyclic colorings is a crucial ingredient of our bidirectional decompression using the augmented matrix.

\subsection{Performance considerations}

In general, we found that the previously available pseudocode did not include all the necessary subtleties for efficient implementation of coloring and decompression algorithms.
Therefore, we have listed some of these missing implementation details, as well as improvements we discovered, in \Cref{sec:implem}.
Notable performance optimizations include:
\begin{itemize}
    \setlength{\itemindent}{-25pt}
    \item Complete avoidance of hash tables for edge index management;
    \item Efficient reconstruction of 2-colored tree structures after acyclic coloring;
    \item Low-memory buckets for the dynamic degree-based orderings.
\end{itemize}

\subsection{Integration with other tools}

\texttt{SMC} integrates seamlessly with two other Julia packages: \texttt{SparseConnectivityTracer.jl}\footnote{\texttt{\url{https://github.com/adrhill/SparseConnectivityTracer.jl}}} and \texttt{DifferentiationInterface.jl}\footnote{\texttt{\url{https://github.com/JuliaDiff/DifferentiationInterface.jl}}}.
The former provides utilities for sparsity pattern detection, a crucial step in the coloring process of Jacobians and Hessians, while the latter offers standardized access to various AD backends for computing JVPs, VJPs and HVPs.
Together, these packages form a cohesive ecosystem for sparse AD in Julia~\citep{hillSparserBetterFaster2025}.

With the upcoming release of Julia version 1.12, it will become possible to distribute small binaries and shared libraries of Julia packages.
Then, we plan to create a C interface for \texttt{SMC}, aiming to improve its interoperability with other established AD systems such as \texttt{ADOL-C} \citep{griewankAlgorithm755ADOLC1996} and \texttt{CasADi} \citep{anderssonCasADiSoftwareFramework2019}.

\section{Numerical experiments}

We performed extensive numerical experiments comparing \texttt{SMC} with \texttt{ColPack} (or \texttt{CP} for short), both in terms of number of colors and in terms of runtime.
Our test instances come from previous works: we reuse the nonsymmetric matrices from \citet[Tables VI and VII]{gebremedhinColPackSoftwareGraph2013}, and the symmetric matrices from \citet[Table 4.2]{gebremedhinWhatColorYour2005}\footnote{except \texttt{mrng1} and \texttt{mrng2}, which we did not find in the \texttt{SuiteSparse} collection}.
For the experiments below, we also include the Linear Programming and Harwell-Boeing matrices studied by \citet[Tables 5.1 and 5.3]{colemanEfficientComputationSparse1998}\footnote{except \texttt{boeing1} and \texttt{boeing2}, which we did not find in the \texttt{SuiteSparse} collection}.
Statistics on the sizes of all the matrices used in our experiments are provided in~\Cref{tab:matrix_table}. 
We only present our main bicoloring results in this section, the rest of the tables can be found in \Cref{sec:benchmarks}.

\begin{table}

\setlength\tabcolsep{0pt}
\centering
\footnotesize
\begin{threeparttable}
\begin{tabular}{@{\extracolsep{2ex}}*{6}{llrrrr}}
\toprule
\textbf{group} & \textbf{name} & \textbf{sym} & \textbf{rows} & \textbf{cols} & \textbf{nnz} \\
\midrule
LPnetlib & lp\_cre\_a & false & 3516 & 7248 & 18168 \\
LPnetlib & lp\_ken\_11 & false & 14694 & 21349 & 49058 \\
LPnetlib & lp\_ken\_13 & false & 28632 & 42659 & 97246 \\
LPnetlib & lp\_maros\_r7 & false & 3136 & 9408 & 144848 \\
LPnetlib & lp\_cre\_d & false & 8926 & 73948 & 246614 \\
LPnetlib & lp\_ken\_18 & false & 105127 & 154699 & 358171 \\
Bai & af23560 & false & 23560 & 23560 & 484256 \\
Shen & e40r0100 & false & 17281 & 17281 & 553562 \\
vanHeukelum & cage11 & false & 39082 & 39082 & 559722 \\
vanHeukelum & cage12 & false & 130228 & 130228 & 2032536 \\
LPnetlib & lp\_standata & false & 359 & 1274 & 3230 \\
LPnetlib & lp\_scagr25 & false & 471 & 671 & 1725 \\
LPnetlib & lp\_scagr7 & false & 129 & 185 & 465 \\
LPnetlib & lp\_stair & false & 356 & 614 & 4003 \\
LPnetlib & lp\_blend & false & 74 & 114 & 522 \\
LPnetlib & lp\_vtp\_base & false & 198 & 346 & 1051 \\
LPnetlib & lp\_agg & false & 488 & 615 & 2862 \\
LPnetlib & lp\_agg2 & false & 516 & 758 & 4740 \\
LPnetlib & lp\_agg3 & false & 516 & 758 & 4756 \\
LPnetlib & lp\_bore3d & false & 233 & 334 & 1448 \\
LPnetlib & lp\_israel & false & 174 & 316 & 2443 \\
LPnetlib & lp\_tuff & false & 333 & 628 & 4561 \\
LPnetlib & lp\_adlittle & false & 56 & 138 & 424 \\
HB & watt\_2 & false & 1856 & 1856 & 11550 \\
HB & can\_256 & true & 256 & 256 & 2916 \\
HB & can\_268 & true & 268 & 268 & 3082 \\
HB & can\_292 & true & 292 & 292 & 2540 \\
HB & can\_634 & true & 634 & 634 & 7228 \\
HB & can\_715 & true & 715 & 715 & 6665 \\
HB & can\_1054 & true & 1054 & 1054 & 12196 \\
HB & can\_1072 & true & 1072 & 1072 & 12444 \\
HB & impcol\_c & false & 137 & 137 & 411 \\
HB & impcol\_d & false & 425 & 425 & 1339 \\
HB & impcol\_e & false & 225 & 225 & 1308 \\
HB & west0067 & false & 67 & 67 & 294 \\
HB & west0381 & false & 381 & 381 & 2157 \\
HB & west0497 & false & 497 & 497 & 1727 \\
HB & gent113 & false & 113 & 113 & 655 \\
HB & arc130 & false & 130 & 130 & 1282 \\
DIMACS10 & 598a & true & 110971 & 110971 & 1483868 \\
DIMACS10 & 144 & true & 144649 & 144649 & 2148786 \\
DIMACS10 & m14b & true & 214765 & 214765 & 3358036 \\
DIMACS10 & auto & true & 448695 & 448695 & 6629222 \\
\bottomrule
\end{tabular}
\end{threeparttable}
\caption{Structure and dimensions of the test matrices.}
\label{tab:matrix_table}
\vspace{-20pt}

\end{table}

\Cref{tab:article_body} compares the variants of bicoloring in \texttt{ColPack} and \texttt{SMC} in terms of the number of colors returned.
It also reports the number of colors obtained by \citet{colemanEfficientComputationSparse1998} with their bicoloring algorithm, referred to as \texttt{CV}, which we did not run ourselves.
The variant of star bicoloring chosen for \texttt{ColPack} is \texttt{ImplicitCoveringStarBicoloring}, which was named as the best overall in \citet{gebremedhinColPackSoftwareGraph2013}.
For the star and acyclic bicoloring in \texttt{SMC}, the post-processing procedure is used.
The boldfaced value indicates the lowest number of colors in each row.

As we can see, our acyclic bicoloring implementation (which enables decompression by substitution and was not available in \texttt{ColPack}) can give rise to fewer colors than implementations of star bicoloring in \texttt{ColPack} and \texttt{SMC}.
Note however that these results are dependent on the chosen ordering: more details are given in \Cref{sec:benchmarks}.
On most instances, the approach of \citet{colemanEfficientComputationSparse1998} was found to yield fewer colors than \texttt{SMC}, yet their method has other drawbacks.
For instance, it relies on partial column- and row-intersection graphs (built from a partition of the matrix $A$ into two sides), which can be denser than the original bipartite graph we use.

\begin{table}

\setlength\tabcolsep{0pt}
\centering
\footnotesize
\begin{threeparttable}
\begin{tabular}{@{\extracolsep{2ex}}*{6}{lccccc}}
\toprule
 & \multicolumn{3}{c}{direct} & \multicolumn{2}{c}{substitution} \\
\cmidrule{2-4}\cmidrule{5-6}
\textbf{name} & CP & SMC & CV~~ & SMC & CV \\
\midrule
lp\_cre\_a & 18 & 16 &  & \textbf{12} &  \\
lp\_ken\_11 & 6 & 5 &  & \textbf{4} &  \\
lp\_ken\_13 & 5 & \textbf{4} &  & \textbf{4} &  \\
lp\_maros\_r7 & 73 & \textbf{72} &  & \textbf{72} &  \\
lp\_cre\_d & 17 & 15 &  & \textbf{13} &  \\
lp\_ken\_18 & 6 & 5 &  & \textbf{4} &  \\
af23560 & 44 & \textbf{32} &  & 45 &  \\
e40r0100 & 88 & \textbf{87} &  & \textbf{87} &  \\
cage11 & 82 & 81 &  & \textbf{79} &  \\
cage12 & 97 & 96 &  & \textbf{93} &  \\
lp\_standata & 11 & 10 & 9 & 8 & \textbf{7} \\
lp\_scagr25 & 10 & 9 & 8 & \textbf{5} & \textbf{5} \\
lp\_scagr7 & 10 & 9 & 8 & \textbf{5} & \textbf{5} \\
lp\_stair & 39 & 38 & 36 & 35 & \textbf{29} \\
lp\_blend & 18 & 17 & 16 & \textbf{14} & \textbf{14} \\
lp\_vtp\_base & 14 & 13 & 12 & \textbf{9} & 10 \\
lp\_agg & 44 & 43 & 19 & 38 & \textbf{13} \\
lp\_agg2 & 44 & 43 & 26 & 38 & \textbf{21} \\
lp\_agg3 & 44 & 43 & 27 & 39 & \textbf{21} \\
lp\_bore3d & 29 & 28 & 28 & 25 & \textbf{24} \\
lp\_israel & 137 & 136 & 61 & 107 & \textbf{49} \\
lp\_tuff & 27 & 25 & 21 & 26 & \textbf{16} \\
lp\_adlittle & 12 & 11 & 11 & 11 & \textbf{10} \\
watt\_2 & 66 & 65 & 20 & \textbf{11} & 12 \\
can\_256 & 84 & 83 & 32 & 27 & \textbf{23} \\
can\_268 & 39 & 38 & 18 & 26 & \textbf{12} \\
can\_292 & 37 & 36 & \textbf{17} & 36 & \textbf{17} \\
can\_634 & 32 & 31 & 28 & 33 & \textbf{21} \\
can\_715 & 106 & 105 & 22 & \textbf{18} & \textbf{18} \\
can\_1054 & 37 & 36 & 31 & 25 & \textbf{23} \\
can\_1072 & 37 & 36 & 32 & 25 & \textbf{24} \\
impcol\_c & 12 & 11 & 6 & 6 & \textbf{4} \\
impcol\_d & 12 & 12 & 6 & 6 & \textbf{5} \\
impcol\_e & 33 & 32 & 21 & 24 & \textbf{14} \\
west0067 & 15 & 14 & 9 & 8 & \textbf{7} \\
west0381 & 75 & 74 & 12 & 73 & \textbf{9} \\
west0497 & 56 & 55 & 22 & \textbf{18} & 19 \\
gent113 & 32 & 31 & 19 & 30 & \textbf{13} \\
arc130 & 125 & 124 & 25 & 125 & \textbf{23} \\
\bottomrule
\end{tabular}
\end{threeparttable}
\caption{Number of colors for bidirectional coloring with natural order, by decompression mode and by package.}
\label{tab:article_body}
\vspace{-20pt}

\end{table}

\section{Conclusion}

We presented new star and acyclic bicoloring algorithms for bidirectional compression-based computation of sparse Jacobian matrices.
They allow the combined use of the forward and reverse modes of AD.
We also presented accompanying efficient decompression procedures.
Our bicoloring and decompression algorithms leverage the technology previously developed for sparse symmetric matrix (Hessian) computation.
Our implementations are delivered via the Julia package \texttt{SparseMatrixColorings.jl}.
This library also includes implementations of coloring and vertex ordering techniques for unidirectional computation of Jacobians (partial distance-2 coloring) and computation of Hessians (star and acyclic coloring).
We showed that \texttt{SMC} exhibits comparable performance with \texttt{ColPack} while offering additional features and reduced code complexity.

We conclude by pointing out a few directions for future work.
\begin{itemize}
\item As observed in our experiments, the number of colors obtained with our new bicoloring technique still compares unfavorably to the algorithm of \citet{colemanEfficientComputationSparse1998}. Understanding the cause of the difference in performance between the two approaches is an important challenge for further research.
\item After the completion of the post-processing step (described in Section 3.1) that neutralizes unneeded colors, the rows and columns of the Jacobian that remain with non-neutral (actual) colors define a vertex cover in the bipartite graph view of the Jacobian. The characterization of this vertex cover in graph-theoretic terms is a worthy subject to investigate in future work.
\item The augmented symmetric matrix for a non-symmetric matrix that we relied on for formulating our star and acyclic bicoloring is known in the literature as the adjoint graph representation of a hypergraph. It is worthwhile to explore this connection for studying hypergraph coloring problems.
\item Our array-based implementation of the recovery routines associated with star and acyclic bicoloring naturally lends itself for parallelization within a SIMD framework. It will be interesting to exploit this to devise efficient GPU implementations.   
\end{itemize}                  

\subsection*{Acknowledgements}

Both first authors contributed equally.

We acknowledge the ideas and feedback contributed by Paul Hovland during the inception of this paper.
We are grateful to Adrian Hill for his implementation of coloring visualization, and for his key role in shaping the sparse AD ecosystem in Julia through joint development efforts.
We sincerely thank Jean-Baptiste Caillau for his assistance in describing the applications in optimal control.

\small
\bibliographystyle{abbrvnat}
\bibliography{bicoloring}
\normalsize

\appendix

\section{Implementation details} \label{sec:implem}

We describe some implementation details for \texttt{SMC}, which improve or extend the algorithms presented in \citet{gebremedhinWhatColorYour2005,gebremedhinNewAcyclicStar2007,gebremedhinEfficientComputationSparse2009,gebremedhinColPackSoftwareGraph2013}.
The detailed pseudocode of every algorithm is given in \Cref{sec:pseudocode}.

\subsection{Graph storage}

\subsubsection{Adjacency graph}

The adjacency graph $\mathcal{G}_a$ of a symmetric matrix $A$ is stored as a Compressed Sparse Column (CSC) encoding with the help of two integer vectors \texttt{colptr} and \texttt{rowval}.
We assume that the row indices in \texttt{rowval} are sorted in increasing order within each column, as is commonly done in sparse matrix libraries to facilitate traversal and efficient access.
Note that by symmetry of $A$, it is also a Compressed Sparse Row (CSR) encoding of $A$.
A minor issue is posed by non-zero diagonal coefficients: they may be present in the original matrix, but they are not explicitly represented in the adjacency graph (which has no self-loops by definition).
Thus, each iteration over neighbors of a vertex $i$ must remember to skip $i$ itself if encountered.
Additionally, for star and acyclic colorings, we need a mapping from edges to unique integer indices.
To achieve this, we construct a vector \texttt{edge\_to\_index} of the same length as the number of nonzeros in $A$, assigning a unique integer $k$ to each edge $(i, j)$ with $i \neq j$.
In particular, this vector ensures that the off-diagonal entries of $\hat{A}$, a sparse CSC matrix defined by \texttt{colptr}, \texttt{rowval}, and \texttt{nzval=edge\_to\_index}, satisfy $\hat{A}_{ij} = \hat{A}_{ji} = k$.
\Cref{alg:edge_to_index} describes how to build this vector from the sparsity pattern of $A$.
This approach is faster for grouped neighbor accesses than a container that stores unique key-value pairs $(i, j) \mapsto k$, such as a dictionary.
Conversely, \texttt{ColPack} implements the mapping between edges and indices using a \texttt{C++} \texttt{map}, which incurs additional overhead due to dynamic memory allocation and tree-based lookups.
For the remainder of the algorithms, we regroup these three vectors into a structure called \texttt{adjacency\_graph} and, for a given vertex \texttt{v}, we use the procedures \texttt{neighbors(adjacency\_graph, v)} and \texttt{indexed\_neighbors(adjacency\_graph, v)} to return the neighbors (with the associated edge indices if needed) other than \texttt{v}.

\subsubsection{Bipartite graph}

For the bipartite graph $\mathcal{G}_b$, dual CSC and CSR encodings of $A$ are necessary for optimal performance.
The CSC version efficiently maps each column $j$ to its row neighbors, and the CSR version efficiently maps each row $i$ to its column neighbors.
At the cost of twice the initial memory space (when $A$ is not structurally symmetric), this allows efficient iteration over neighbors for both parts of the vertex space.
To save memory, both CSC and CSR structures are manipulated without storing the nonzero values of $A$.
\Cref{alg:transpose_sparsity_pattern} details the conversion process between the CSC and CSR formats of $A$.
We assume that row indices are sorted in increasing order within each column for the conversion.

\subsection{Colorings}

\subsubsection{Partial distance-2 coloring}

We implement Algorithm 3.2 of \citet{gebremedhinWhatColorYour2005} without any modification.

\subsubsection{Star coloring}

We implement Algorithm 4.1 of \citet{gebremedhinNewAcyclicStar2007} without using hash tables, as detailed in~\Cref{alg:star_coloring}.
The original algorithm did not assign a hub to trivial stars. However, the corresponding decompression routine requires that every star have a designated hub.
To address this, we choose an arbitrary hub for trivial stars but encode it as a negative value $-h$.
This convention not only makes trivial stars easily distinguishable from non-trivial ones during star coloring and post-processing but also simplifies the decompression process.

\subsubsection{Acyclic coloring}

We implement Algorithm 3.1 of \citet{gebremedhinNewAcyclicStar2007} without using hash tables, as detailed in~\Cref{alg:acyclic_coloring}.
To efficiently build the 2-colored trees, we use a disjoint-set structure that requires two integer vectors \texttt{parents} and \texttt{ranks}, as well as an integer \texttt{nt} to keep track of the current number of trees.
In the rest of the paper, this data structure is referred to as \texttt{forest}.
Three essential routines are implemented on this structure: \texttt{create\_forest}, \texttt{find\_root} and \texttt{root\_union}.
They are detailed in~\Cref{alg:forest}.
Note that the edge index used to represent the root in the structure \texttt{forest} does not correspond to the actual tree root vertex.
Instead, the selected ``root'' is maintained to minimize the tree’s height (i.e., it is chosen following a union-by-rank strategy to optimize path compression).

\subsubsection{Star and acyclic bicoloring}

To perform star and acyclic bicoloring, we build the adjacency graph of the augmented symmetric matrix $H$~\eqref{eq:augmented_matrix}, derived from a rectangular matrix $J$.
This can be done by efficiently combining \Cref{alg:edge_to_index} and \Cref{alg:transpose_sparsity_pattern}.
To compute the seeds for the products $Ju$ and $J^{\top\!}v$, we employ the vectors \texttt{row\_colors} and \texttt{column\_colors} produced by \Cref{alg:remap_colors}.
Additionally, we use \texttt{sym\_to\_row} and \texttt{sym\_to\_col} (also generated by \Cref{alg:remap_colors}) along with the relation~\eqref{eq:Bs_Br_Bc} for decompression.

\subsubsection{Post-processing}

Detecting the trivial structures among the stars and trees is simple because we only need to check whether the number of vertices in each structure is exactly two (or one edge).
However, in acyclic coloring, we need to distinguish the stars from the normal trees.
To do that, we determine a boolean vector \texttt{is\_star} that indicates for each tree whether it is a star (potentially trivial) or not.
We can compute this concurrently with the reverse BFS orders of the two-colored trees (\Cref{alg:reverse_bfs}).

Note that we can reuse the vector \texttt{forbidden\_colors} needed in both star and acyclic colorings (Algorithms 3.1 and 4.1 of \citet{gebremedhinNewAcyclicStar2007}) for the vector \texttt{offsets} during post-processing (\Cref{alg:postprocessing}).

\subsection{Decompression}

In \texttt{SMC}, we optimize the decompression of all our colorings for sparse CSC matrices with a generic fallback for matrices in other formats.
Except for row and column coloring, where the decompression is straightforward, the decompression of symmetric colorings and derived bicolorings uses the algorithms \textit{DirectRecover2} and \textit{IndirectRecover} from \citet{gebremedhinEfficientComputationSparse2009}, combined with new additional strategies.

\subsubsection{Direct decompression}

For direct decompression, we found it more efficient to use linear indexing and precompute the positions of nonzeros in the compressed matrix $B$.
We store them in a dedicated vector \texttt{compressed\_indices}.
$$
    \texttt{nzval[k] = B[compressed\_indices[k]]} \quad \text{for} \quad \texttt{k} = 1, \dots, \texttt{nnzA}.
$$
For example, after column coloring, a coefficient $A_{ij}$ is recovered from $B_{ic}$ where \texttt{c = col\_colors[j]} at the linear index $(c-1) \times m + i$ when $B$ is column-ordered and 1-based indexed.
Although the recovery is simple for row or column coloring, it becomes more complex and less parallel for star coloring and bicoloring.
We need to use the stars to determine which vertex's color of a edge $(i,j)$ is needed to recover a coefficient $A_{ij}$.
Note that an additional integer vector \texttt{A\_indices} is needed for star bicoloring to specify which nonzeros are recovered from $B_r$ and $B_c$:
\begin{align*}
    \texttt{nzval[A\_indices[k]]} & = \texttt{$B_c$[compressed\_indices[k]]}, \quad \text{for} \quad \texttt{k} = 1, \dots, \ell,               \\
    \texttt{nzval[A\_indices[k]]} & = \texttt{$B_r$[compressed\_indices[k]]}, \quad \text{for} \quad \texttt{k} = \ell+1, \dots, \texttt{nnzA}.
\end{align*}
The variable $\ell$ denotes the number of nonzeros of $A$ recovered from $B_c$.

\subsubsection{Decompression by substitution}

For acyclic coloring and bicoloring, the decompression step relies on traversing the two-colored trees.
However, existing work does not explain how to reconstruct the structure of these trees from their disjoint-set representation.
While the coloring phase identifies which edges belong to the same tree, it does not reveal how each tree can be traversed from its leaves to its root, following a reverse breadth-first search (BFS) order.
Our goal is to clarify how this ordering can be computed and extracted from the structure \texttt{forest}.
Once determined, the edges within each tree are sorted according to the order in which they must be processed.

In this ordering, each edge is oriented so that its first vertex is a leaf (after pruning previously processed edges), and its second vertex is an internal node.
This convention enables the direct application of \textit{IndirectRecover} from \citet{gebremedhinEfficientComputationSparse2009}.

The reconstruction of the trees is described in~\Cref{alg:structure_trees}, and the computation of the reverse BFS order is presented in~\Cref{alg:reverse_bfs}.
We use the fact that a tree with \texttt{ne} edges has exactly \texttt{ne+1} vertices in these routines.
Note that the vector \texttt{ranks} from the structure \texttt{forest} can be repurposed to store \texttt{root\_to\_tree} in~\Cref{alg:structure_trees}.
Likewise, the array \texttt{first\_visit\_to\_tree} used during acyclic coloring can be reused to store \texttt{reverse\_bfs\_orders} in~\Cref{alg:reverse_bfs}.

With this approach, the decompression of each tree can be performed in parallel.  
However, it offers less parallelism than direct decompression, since parallelism occurs across the number of trees rather than across the number of nonzeros.  
Moreover, the number of edges per tree may vary significantly, which may require load balancing to efficiently exploit this parallelism.

For this reason, we emphasize that \textit{IndirectRecover} from \citet{gebremedhinEfficientComputationSparse2009} is equivalent to solving a sparse triangular system in which all coefficients are either 0 or 1, the unknowns correspond to the nonzeros, and the right-hand side is a subset of the entries of the matrix $B$.
An alternative is to explicitly compute the inverses of all the sparse triangular matrices needed.
These inverses remain triangular and contain only -1 (off-diagonal), 0, and 1 (diagonal) as entries.
All of these inverse matrices can then be assembled into a single large sparse matrix, and the full vector of \texttt{nzval} can be recovered through a single sparse matrix-vector multiplication with \texttt{B[compressed\_indices]}.
This approach might require significantly more storage, but it is likely the most efficient way to exploit SIMD architectures such as GPUs.

\subsection{Vertex orderings} \label{sec:orders}

Our implementation of vertex orderings follows Section 5.1 of \citet{gebremedhinColPackSoftwareGraph2013}.
The authors describe three different orderings using a common framework, that of dynamic degree-based orderings.
Each vertex in the graph gets assigned a dynamic degree, which evolves as the ordering progresses and the permutation $\pi$ of vertices gets built.
More precisely, the back degree of a vertex $v$ is the number of neighbors that appear before $v$ in the permutation $\pi$, while the forward degree of $v$ is the number of neighbors that appear after $v$.
As the permutation $\pi = \{v_1, \dots, v_n\}$ is grown iteratively, either from $v_1$ to $v_n$ or from $v_n$ to $v_1$, the dynamic degrees (back or forward) evolve.
Since these degrees are used to select which vertex will be added to the permutation next, we need an efficient way to keep track of them.

\citet[Figure 2]{gebremedhinColPackSoftwareGraph2013} propose a bucket sort mechanism to achieve this goal.
They use one variable-size bucket per possible value of the dynamic degree.
Whenever the dynamic degree of a vertex changes, it gets moved one bucket up or one bucket down, and appended to the end of the destination bucket.
This is rather straightforward to implement, but requires one stack-like data structure per bucket, whose size will increase or decrease and cannot be predicted ahead of time.

To improve memory locality and remove the need for variable-sized containers, we suggest an alternative where all buckets are stored inside the same vector, and the limits of each bucket (an interval of integers) are kept updated separately.
The key is to allow modification of a bucket either from the start or from the end.
That way, when a vertex $v$ moves up from bucket $d$ to bucket $d+1$, bucket $d$ shrinks from the right and bucket $d+1$ grows to the left.
When a vertex $v$ moves down from bucket $d$ to bucket $d-1$, bucket $d$ shrinks from the left and bucket $d-1$ grows to the right.
In both cases, the total amount of storage necessary remains constant and predictable.

The two algorithms are detailed and compared in \Cref{alg:increase_degree_stacks} and \Cref{alg:increase_degree_vector}, on the special case of a vertex whose dynamic degree needs to increase by one.

\section{Benchmarks} \label{sec:benchmarks}

\subsection{Setup}

The benchmarks are run on a 14-inch MacBook Pro M3 (2023) with 36 GB of RAM, using Julia v1.11.5 in single-threaded mode.
\texttt{SMC} is called through its native Julia interface, which accepts CSC matrices, while \texttt{ColPack} is called through its command line interface, which accepts \texttt{.mtx} matrix files.
We use a compiled version of \texttt{ColPack} without OpenMP (through \texttt{ColPack.jl} v0.5.0), along with \texttt{SparseMatrixColorings.jl} v0.4.19.

\subsection{Test matrices}

The test matrices are retrieved from the \texttt{SuiteSparse} matrix collection \citep{davisUniversityFloridaSparse2011,kolodziejSuiteSparseMatrixCollection2019}.
We describe their main features in \Cref{tab:matrix_table}.

\subsection{Number of colors}

We first investigate whether or not our \texttt{SMC} implementation is consistent with the results returned by \texttt{ColPack}.
This is done for each of the following orders: smallest last (SL), incidence degree (ID), largest first (LF), dynamic largest first (DLF) and natural (N).
Boldfaced numbers indicate the minimum number of colors for a given instance across all orders.

We recall that some implementation differences between \texttt{SMC} and \texttt{ColPack} can affect the number of colors:
\begin{itemize}
    \item \texttt{SMC} has a completely different algorithm for star bicoloring (\Cref{sec:bitosym});
    \item The natural order for bicoloring differs: \texttt{SMC} colors columns before rows, whereas \texttt{ColPack} does the opposite;
    \item \texttt{SMC} calls a post-processing routine after star and acyclic coloring algorithms (\Cref{sec:postprocessing}), which can reduce the number of colors if some diagonal coefficients are zero.
\end{itemize}
To minimize discrepancies, we always imitate \texttt{ColPack}'s stack-based ordering algorithm instead of leveraging our own vector-based improvement (\Cref{sec:orders}).

\subsubsection{Unidirectional and symmetric colorings}

Results for row coloring are presented in \Cref{tab:results_row_direct}, for column coloring in \Cref{tab:results_column_direct}, for star coloring in \Cref{tab:results_symmetric_direct}, for acyclic coloring in \Cref{tab:results_symmetric_substitution}.
In these four tables, underlined numbers indicate a discrepancy between \texttt{SMC} and \texttt{ColPack}.
Each table contains at most one instance where the two libraries differ (usually by a couple of colors), which suggests that both implementations give very similar results.
\begin{table}[htb!]
    
\setlength\tabcolsep{0pt}
\centering
\footnotesize
\begin{threeparttable}
\begin{tabular}{@{\extracolsep{2ex}}*{11}{lcccccccccc}}
\toprule
 & \multicolumn{2}{c}{SL} & \multicolumn{2}{c}{ID} & \multicolumn{2}{c}{LF} & \multicolumn{2}{c}{DLF} & \multicolumn{2}{c}{N} \\
\cmidrule{2-3}\cmidrule{4-5}\cmidrule{6-7}\cmidrule{8-9}\cmidrule{10-11}
\textbf{name} & CP & SMC & CP & SMC & CP & SMC & CP & SMC & CP & SMC \\
\midrule
lp\_cre\_a & \textbf{14} & \textbf{14} & \textbf{14} & \textbf{14} & \textbf{14} & \textbf{14} & \textbf{14} & \textbf{14} & 16 & 16 \\
lp\_ken\_11 & \textbf{4} & \textbf{4} & \textbf{4} & \textbf{4} & \textbf{4} & \textbf{4} & 5 & 5 & 5 & 5 \\
lp\_ken\_13 & \textbf{4} & \textbf{4} & 5 & 5 & 5 & 5 & 5 & 5 & \textbf{4} & \textbf{4} \\
lp\_maros\_r7 & 80 & 80 & 88 & 88 & 100 & 100 & 113 & 113 & \textbf{72} & \textbf{72} \\
lp\_cre\_d & 15 & 15 & 15 & 15 & \textbf{13} & \textbf{13} & 14 & 14 & 15 & 15 \\
lp\_ken\_18 & \textbf{4} & \textbf{4} & 5 & 5 & 5 & 5 & 5 & 5 & 5 & 5 \\
af23560 & 42 & 42 & 43 & \underline{42} & 43 & \underline{36} & 60 & \underline{59} & 43 & \textbf{\underline{32}} \\
e40r0100 & \textbf{70} & \textbf{70} & 71 & 71 & 85 & 85 & 86 & 86 & 87 & 87 \\
cage11 & \textbf{64} & \textbf{64} & 67 & 67 & 67 & 67 & 70 & 70 & 81 & 81 \\
cage12 & \textbf{67} & \textbf{67} & 72 & 72 & 73 & 73 & 79 & 79 & 96 & 96 \\
\bottomrule
\end{tabular}
\end{threeparttable}
\caption{Number of colors for row coloring.}
\label{tab:results_row_direct}
\vspace{-20pt}

\end{table}
\begin{table}[htb!]
    
\setlength\tabcolsep{0pt}
\centering
\footnotesize
\begin{threeparttable}
\begin{tabular}{@{\extracolsep{2ex}}*{11}{lcccccccccc}}
\toprule
 & \multicolumn{2}{c}{SL} & \multicolumn{2}{c}{ID} & \multicolumn{2}{c}{LF} & \multicolumn{2}{c}{DLF} & \multicolumn{2}{c}{N} \\
\cmidrule{2-3}\cmidrule{4-5}\cmidrule{6-7}\cmidrule{8-9}\cmidrule{10-11}
\textbf{name} & CP & SMC & CP & SMC & CP & SMC & CP & SMC & CP & SMC \\
\midrule
lp\_cre\_a & \textbf{360} & \textbf{360} & \textbf{360} & \textbf{360} & \textbf{360} & \textbf{360} & \textbf{360} & \textbf{360} & \textbf{360} & \textbf{360} \\
lp\_ken\_11 & 125 & 125 & 124 & 124 & 128 & 128 & \textbf{122} & \textbf{122} & 130 & 130 \\
lp\_ken\_13 & 171 & 171 & 171 & 171 & 174 & 174 & \textbf{170} & \textbf{170} & 176 & 176 \\
lp\_maros\_r7 & 83 & 83 & 90 & 90 & \textbf{70} & \textbf{70} & 114 & 114 & 74 & 74 \\
lp\_cre\_d & \textbf{808} & \textbf{808} & \textbf{808} & \textbf{808} & \textbf{808} & \textbf{808} & \textbf{808} & \textbf{808} & 813 & 813 \\
lp\_ken\_18 & \textbf{325} & \textbf{325} & 326 & 326 & 328 & 328 & \textbf{325} & \textbf{325} & 330 & 330 \\
af23560 & 42 & 42 & 42 & 42 & 43 & \underline{44} & 59 & 59 & 44 & \textbf{\underline{32}} \\
e40r0100 & \textbf{70} & \textbf{70} & 71 & 71 & 87 & 87 & 85 & 85 & 95 & 95 \\
cage11 & \textbf{64} & \textbf{64} & 67 & 67 & 67 & 67 & 70 & 70 & 81 & 81 \\
cage12 & \textbf{67} & \textbf{67} & 72 & 72 & 73 & 73 & 79 & 79 & 96 & 96 \\
\bottomrule
\end{tabular}
\end{threeparttable}
\caption{Number of colors for column coloring.}
\label{tab:results_column_direct}
\vspace{-20pt}

\end{table}
\begin{table}[htb!]
    
\setlength\tabcolsep{0pt}
\centering
\footnotesize
\begin{threeparttable}
\begin{tabular}{@{\extracolsep{2ex}}*{11}{lcccccccccc}}
\toprule
 & \multicolumn{2}{c}{SL} & \multicolumn{2}{c}{ID} & \multicolumn{2}{c}{LF} & \multicolumn{2}{c}{DLF} & \multicolumn{2}{c}{N} \\
\cmidrule{2-3}\cmidrule{4-5}\cmidrule{6-7}\cmidrule{8-9}\cmidrule{10-11}
\textbf{name} & CP & SMC & CP & SMC & CP & SMC & CP & SMC & CP & SMC \\
\midrule
598a & \textbf{23} & \textbf{23} & \textbf{23} & \textbf{23} & 27 & 27 & 28 & 28 & 28 & 28 \\
144 & \textbf{24} & \textbf{24} & 25 & 25 & 28 & 28 & 30 & 30 & 29 & 29 \\
m14b & 28 & 28 & \textbf{27} & \textbf{27} & 36 & 36 & 38 & 38 & 32 & 32 \\
auto & \textbf{27} & \textbf{27} & 28 & \underline{29} & 36 & 36 & 34 & 34 & 32 & 32 \\
\bottomrule
\end{tabular}
\end{threeparttable}
\caption{Number of colors for star coloring.}
\label{tab:results_symmetric_direct}
\vspace{-20pt}

\end{table}
\begin{table}[htb!]
    
\setlength\tabcolsep{0pt}
\centering
\footnotesize
\begin{threeparttable}
\begin{tabular}{@{\extracolsep{2ex}}*{11}{lcccccccccc}}
\toprule
 & \multicolumn{2}{c}{SL} & \multicolumn{2}{c}{ID} & \multicolumn{2}{c}{LF} & \multicolumn{2}{c}{DLF} & \multicolumn{2}{c}{N} \\
\cmidrule{2-3}\cmidrule{4-5}\cmidrule{6-7}\cmidrule{8-9}\cmidrule{10-11}
\textbf{name} & CP & SMC & CP & SMC & CP & SMC & CP & SMC & CP & SMC \\
\midrule
598a & 12 & 12 & 12 & \textbf{\underline{11}} & 12 & 12 & 12 & 12 & 13 & 13 \\
144 & \textbf{13} & \textbf{13} & \textbf{13} & \textbf{13} & \textbf{13} & \textbf{13} & 14 & 14 & 14 & 14 \\
m14b & \textbf{14} & \textbf{14} & \textbf{14} & \textbf{14} & 15 & 15 & 16 & 16 & 18 & 18 \\
auto & \textbf{14} & \textbf{14} & 15 & 15 & 15 & 15 & 16 & 16 & 16 & 16 \\
\bottomrule
\end{tabular}
\end{threeparttable}
\caption{Number of colors for acyclic coloring.}
\label{tab:results_symmetric_substitution}
\vspace{-20pt}

\end{table}

\subsubsection{Bidirectional colorings}

Results for star bicoloring are presented in \Cref{tab:results_bidirectional_direct}, and for acyclic bicoloring in \Cref{tab:results_bidirectional_substitution}.
Here, we observe that the lowest number of colors is nearly always provided by an \texttt{SMC} implementation, and that acyclic bicoloring leads to fewer colors than star bicoloring, as expected.
In both star and acyclic bicoloring, the number of colors produced by the \texttt{SMC} implementation is heavily dependent on the ordering, and some orderings perform very poorly.
We conjecture that this may be due to the orderings being applied to the augmented symmetric matrix, whereas they were designed for a different graph representation. Designing orderings suitable for the \texttt{SMC} approach is an interesting direction for future research.
\begin{table}[htb!]
    
\setlength\tabcolsep{0pt}
\centering
\footnotesize
\begin{threeparttable}
\begin{tabular}{@{\extracolsep{2ex}}*{11}{lcccccccccc}}
\toprule
 & \multicolumn{2}{c}{SL} & \multicolumn{2}{c}{ID} & \multicolumn{2}{c}{LF} & \multicolumn{2}{c}{DLF} & \multicolumn{2}{c}{N} \\
\cmidrule{2-3}\cmidrule{4-5}\cmidrule{6-7}\cmidrule{8-9}\cmidrule{10-11}
\textbf{name} & CP & SMC & CP & SMC & CP & SMC & CP & SMC & CP & SMC \\
\midrule
lp\_cre\_a & 27 & 31 & 24 & 370 & 23 & 364 & 21 & 365 & 18 & \textbf{16} \\
lp\_ken\_11 & 128 & 68 & 129 & 128 & 8 & 134 & \textbf{5} & 133 & 6 & \textbf{5} \\
lp\_ken\_13 & 176 & 166 & 177 & 169 & 8 & 180 & 6 & 182 & 5 & \textbf{4} \\
lp\_maros\_r7 & 131 & 80 & 153 & 91 & 130 & 134 & \textbf{71} & 116 & 73 & 72 \\
lp\_cre\_d & 38 & 275 & 25 & 815 & 23 & 811 & 17 & 810 & 17 & \textbf{15} \\
lp\_ken\_18 & 331 & 241 & 331 & 332 & 8 & 334 & 6 & 337 & 6 & \textbf{5} \\
af23560 & 75 & 53 & 74 & 77 & 67 & 71 & 51 & 105 & 44 & \textbf{32} \\
e40r0100 & 99 & 155 & 98 & \textbf{72} & 93 & 83 & 89 & 147 & 88 & 87 \\
cage11 & 96 & 76 & 105 & 81 & 94 & 87 & \textbf{69} & 158 & 82 & 81 \\
cage12 & 102 & 106 & 114 & 90 & 104 & 97 & \textbf{74} & 121 & 97 & 96 \\
\bottomrule
\end{tabular}
\end{threeparttable}
\caption{Number of colors for star bicoloring.}
\label{tab:results_bidirectional_direct}
\vspace{-20pt}

\end{table}
\begin{table}[htb!]
    
\setlength\tabcolsep{0pt}
\centering
\footnotesize
\begin{threeparttable}
\begin{tabular}{@{\extracolsep{2ex}}*{6}{lccccc}}
\toprule
 & SL & ID & LF & DLF & N \\
\textbf{name} & SMC & SMC & SMC & SMC & SMC \\
\midrule
lp\_cre\_a & 43 & 118 & 116 & 117 & \textbf{12} \\
lp\_ken\_11 & 10 & 11 & 12 & 12 & \textbf{4} \\
lp\_ken\_13 & 10 & 10 & 11 & 10 & \textbf{4} \\
lp\_maros\_r7 & 159 & 97 & 115 & 115 & \textbf{72} \\
lp\_cre\_d & 87 & 152 & 151 & 148 & \textbf{13} \\
lp\_ken\_18 & 10 & 10 & 9 & 10 & \textbf{4} \\
af23560 & 60 & 61 & 57 & 87 & \textbf{45} \\
e40r0100 & 157 & \textbf{72} & 83 & 148 & 87 \\
cage11 & 46 & \textbf{45} & 52 & 125 & 79 \\
cage12 & 47 & \textbf{45} & 55 & 99 & 93 \\
\bottomrule
\end{tabular}
\end{threeparttable}
\caption{Number of colors for acyclic bicoloring (\texttt{SMC} only).}
\label{tab:results_bidirectional_substitution}
\vspace{-20pt}

\end{table}

\subsection{Timings}

We have established that the \texttt{SMC} and \texttt{ColPack} algorithms are, at minimum, directly comparable.
The next question is one of performance: can a Julia implementation compete with the speed of an optimized \texttt{C++} library?
To answer it, we run the function of interest five times and record the minimum runtime for each package.
Then, we display the ratio of these minimum runtimes, \texttt{SMC} over \texttt{ColPack}.
Ratios under $1$ indicates that \texttt{SMC} is faster, while ratios above $1$ show that \texttt{ColPack} stays ahead.
We do not measure the time it takes to process the matrix and construct the graph (which is format-dependent): we only measure the ordering time and the subsequent coloring time once the graph is available.

\subsubsection{Coloring times}

We start with the results for coloring times in the case of the natural order, presented in \Cref{tab:color_times_nonsym} and \Cref{tab:color_times_sym}.
The unidirectional colorings in \texttt{SMC} are generally as fast as \texttt{ColPack}'s, 
but our symmetric and bidirectional colorings can be up to four times faster.
The slight overhead in the unidirectional routines stems from Julia’s bounds checking on array accesses, an operation not present in \texttt{C++}, which may modestly constrain compiler optimizations.
Note that for acyclic bicoloring, we compute the runtime ratio compared to \texttt{ColPack}'s star bicoloring, for lack of an equivalent algorithm.
\begin{table}[htb!]
    \begin{minipage}{0.55\textwidth}
        
\setlength\tabcolsep{0pt}
\centering
\footnotesize
\begin{threeparttable}
\begin{tabular}{@{\extracolsep{2ex}}*{5}{lcccc}}
\toprule
 & column & row & \multicolumn{2}{c}{bidirectional$^\star$} \\
\cmidrule{4-5}
\textbf{name} &  &  & star & acyclic \\
\midrule
lp\_cre\_a & 1.2 & 0.808 & 0.27 & 0.563 \\
lp\_ken\_11 & 1.11 & 1.08 & 0.202 & 0.424 \\
lp\_ken\_13 & 1.17 & 1.08 & 0.22 & 0.405 \\
lp\_maros\_r7 & 0.899 & 1.13 & 0.568 & 1.22 \\
lp\_cre\_d & 1.37 & 1.16 & 0.227 & 0.497 \\
lp\_ken\_18 & 1.22 & 1.22 & 0.201 & 0.442 \\
af23560 & 1.35 & 1.43 & 0.377 & 0.973 \\
e40r0100 & 1.17 & 1.31 & 0.615 & 1.22 \\
cage11 & 1.23 & 1.28 & 0.414 & 1.2 \\
cage12 & 1.28 & 1.33 & 0.393 & 1.28 \\
\bottomrule
\end{tabular}
\end{threeparttable}
\caption{Coloring time ratios\\(\texttt{SMC}/\texttt{ColPack}) for row coloring, column coloring and bicoloring with natural order.\\ {\footnotesize$^\star$ both compared to star bicoloring in \texttt{ColPack}}}
\label{tab:color_times_nonsym}
\vspace{-20pt}

    \end{minipage}
    \hfill
    \begin{minipage}{0.38\textwidth}
        
\setlength\tabcolsep{0pt}
\centering
\footnotesize
\begin{threeparttable}
\begin{tabular}{@{\extracolsep{2ex}}*{3}{lcc}}
\toprule
 & \multicolumn{2}{c}{unidirectional} \\
\cmidrule{2-3}
\textbf{name} & star & acyclic \\
\midrule
598a & 0.271 & 0.449 \\
144 & 0.267 & 0.455 \\
m14b & 0.277 & 0.449 \\
auto & 0.304 & 0.489 \\
\bottomrule
\end{tabular}
\end{threeparttable}
\caption{Coloring time ratios\\(\texttt{SMC}/\texttt{ColPack}) for star and acyclic coloring with natural order.}
\label{tab:color_times_sym}
\vspace{-20pt}

    \end{minipage}
\end{table}

\subsubsection{Ordering times}

We present results for the ordering subroutines in \Cref{tab:order_times_row}, \Cref{tab:order_times_column}, and \Cref{tab:order_times_symmetric}.
As we can see, \texttt{SMC} is usually within a factor of $2$ of \texttt{ColPack}'s ordering performance.
Note that for bicoloring, \texttt{SMC} runs the ordering on the augmented adjacency graph, while \texttt{ColPack} runs it on the initial bipartite graph, which makes the timings harder to compare meaningfully.
\begin{table}[htb!]
    \begin{minipage}{0.45\textwidth}
        
\setlength\tabcolsep{0pt}
\centering
\footnotesize
\begin{threeparttable}
\begin{tabular}{@{\extracolsep{2ex}}*{5}{lcccc}}
\toprule
\textbf{name} & SL & ID & LF & DLF \\
\midrule
lp\_cre\_a & 0.83 & 1.04 & 0.84 & 1.06 \\
lp\_ken\_11 & 1.18 & 1.49 & 1.57 & 1.35 \\
lp\_ken\_13 & 1.16 & 1.37 & 1.63 & 1.39 \\
lp\_maros\_r7 & 1.3 & 1.34 & 1.4 & 0.91 \\
lp\_cre\_d & 1.11 & 1.18 & 1.04 & 1.24 \\
lp\_ken\_18 & 1.22 & 1.47 & 1.97 & 1.37 \\
af23560 & 1.46 & 1.73 & 1.89 & 1.05 \\
e40r0100 & 0.87 & 0.97 & 1.21 & 0.86 \\
cage11 & 1.17 & 1.19 & 1.2 & 1.22 \\
cage12 & 1.24 & 1.23 & 1.29 & 1.31 \\
\bottomrule
\end{tabular}
\end{threeparttable}
\caption{Ordering time ratios \\(\texttt{SMC}/\texttt{ColPack}) for row coloring.}
\label{tab:order_times_row}
\vspace{-20pt}

    \end{minipage}
    \hfill
    \begin{minipage}{0.45\textwidth}
        
\setlength\tabcolsep{0pt}
\centering
\footnotesize
\begin{threeparttable}
\begin{tabular}{@{\extracolsep{2ex}}*{5}{lcccc}}
\toprule
\textbf{name} & SL & ID & LF & DLF \\
\midrule
lp\_cre\_a & 1.36 & 0.91 & 1.86 & 1.42 \\
lp\_ken\_11 & 1.82 & 1.29 & 2.41 & 1.46 \\
lp\_ken\_13 & 1.57 & 1.28 & 2.48 & 1.58 \\
lp\_maros\_r7 & 1.38 & 1.22 & 1.36 & 0.94 \\
lp\_cre\_d & 1.64 & 1.39 & 1.58 & 1.58 \\
lp\_ken\_18 & 1.78 & 1.41 & 2.94 & 1.71 \\
af23560 & 1.5 & 1.7 & 1.7 & 1 \\
e40r0100 & 0.95 & 1 & 1.28 & 0.95 \\
cage11 & 1.18 & 1.08 & 1.19 & 1.24 \\
cage12 & 1.24 & 1.12 & 1.26 & 1.31 \\
\bottomrule
\end{tabular}
\end{threeparttable}
\caption{Ordering time ratios\\(\texttt{SMC}/\texttt{ColPack}) for column coloring.}
\label{tab:order_times_column}
\vspace{-20pt}

    \end{minipage}

    \vspace{20pt}
    
\setlength\tabcolsep{0pt}
\centering
\footnotesize
\begin{threeparttable}
\begin{tabular}{@{\extracolsep{2ex}}*{5}{lcccc}}
\toprule
\textbf{name} & SL & ID & LF & DLF \\
\midrule
598a & 1.71 & 1.68 & 4.96 & 1.59 \\
144 & 1.76 & 1.74 & 5.81 & 1.9 \\
m14b & 1.84 & 1.67 & 5.89 & 1.97 \\
auto & 1.94 & 1.83 & 7.41 & 1.73 \\
\bottomrule
\end{tabular}
\end{threeparttable}
\caption{Ordering time ratios\\(\texttt{SMC}/\texttt{ColPack}) for star coloring.}
\label{tab:order_times_symmetric}
\vspace{-20pt}

\end{table}

\clearpage

\section{Pseudocode} \label{sec:pseudocode}

In this appendix, we present the detail of the algorithms we adapted or designed for ordering, coloring and decompression. 

\begin{algorithm}[ht!]
  \footnotesize
  \begin{framed}
    \footnotesize
\textbf{Input:}
\begin{itemize}
  \setlength{\itemsep}{0pt}
  \setlength{\parskip}{0pt}
  \setlength{\itemindent}{-12pt}
  \item An integer \texttt{n} corresponding to the dimension of a symmetric matrix $A$.
  \item An integer \texttt{nnzA} corresponding to the number of non-zeros of $A$.
  \item An integer vector \texttt{colptr} of length \texttt{n+1}.
  \item An integer vector \texttt{rowval} of length \texttt{nnzA}.
  \item An integer vector \texttt{offsets} of length \texttt{n}.
\end{itemize}

\textbf{Output:}
\begin{itemize}
  \setlength{\itemsep}{0pt}
  \setlength{\parskip}{0pt}
  \setlength{\itemindent}{-12pt}
  \item An integer vector \texttt{edge\_to\_index} of length \texttt{nnzA}.
\end{itemize}

\begin{algorithmic}[1]
  \State Initialize \texttt{offsets[1\ldots n]} $\gets$ \texttt{0}
  \State Set \texttt{k} $\gets$ \texttt{0}
  \For{\texttt{j} $\gets$ \texttt{1} \textbf{to} \texttt{n}}
  \For{\texttt{p} $\gets$ \texttt{colptr[j]} \textbf{to} \texttt{colptr[j+1]-1}}
  \State \texttt{i} $\gets \texttt{rowval[p]}$
  \If{\texttt{i > j}}
  \State \texttt{k} $\gets$ \texttt{k + 1}
  \State \texttt{edge\_to\_index[p]} $\gets$ \texttt{k}
  \State \texttt{q} $\gets$ \texttt{colptr[i] + offsets[i]}
  \State \texttt{edge\_to\_index[q]} $\gets$ \texttt{k}
  \State \texttt{offsets[i]} $\gets$ \texttt{offsets[i] + 1}
  \EndIf
  \EndFor
  \EndFor
  \State \Return \texttt{edge\_to\_index}
\end{algorithmic}

  \end{framed}
    \caption{Procedure to build the vector \texttt{edge\_to\_index} from a sparse symmetric matrix $A$ stored in CSC format, where the entries in \texttt{rowval} for each column must be sorted in ascending order.
    The arrays use 1-based indexing.}
  \label{alg:edge_to_index}
\end{algorithm}

\begin{algorithm}[ht!]
  \footnotesize
  \begin{framed}
    \textbf{Input:}
\begin{itemize}
  \setlength{\itemsep}{0pt}
  \setlength{\parskip}{0pt}
  \setlength{\itemindent}{-12pt}
  \item An integer \texttt{m} corresponding to the number of rows of $A$.
  \item An integer \texttt{n} corresponding to the number of columns of $A$.
  \item An integer \texttt{nnzA} corresponding to the number of non-zeros of $A$.
  \item An integer vector \texttt{colptr} of length \texttt{n+1}.
  \item An integer vector \texttt{rowval} of length \texttt{nnzA}.
\end{itemize}

\textbf{Output:}
\begin{itemize}
  \setlength{\itemsep}{0pt}
  \setlength{\parskip}{0pt}
  \setlength{\itemindent}{-12pt}
  \item An integer vector \texttt{rowptr} of length \texttt{m+1}.
  \item An integer vector \texttt{colval} of length \texttt{nnzA}.
\end{itemize}

\begin{algorithmic}[1]
  \State Initialize \texttt{rowptr[1\ldots m+1]} $\gets$ \texttt{0}
  \For{$k \gets$ \texttt{1} \textbf{to} \texttt{nnzA}}
  \State \texttt{i} $\gets$ \texttt{rowval[k]}
  \State \texttt{rowptr[i]} $\gets$ \texttt{rowptr[i] + 1}
  \EndFor
  \State Set \texttt{counter} $\gets$ \texttt{1}
  \For{$i \gets$ \texttt{1} \textbf{to} \texttt{m}}
  \State \texttt{nnz\_row} $\gets$ \texttt{rowptr[i]}
  \State \texttt{rowptr[i]} $\gets$ \texttt{counter}
  \State \texttt{counter} $\gets$ \texttt{counter + nnz\_row}
  \EndFor
  \State \texttt{rowptr[m+1]} $\gets$ \texttt{counter}

  \For{$j \gets$ \texttt{1} \textbf{to} \texttt{n}}
  \For{$p \gets$ \texttt{colptr[j]} \textbf{to} \texttt{colptr[j+1]-1}}
  \State \texttt{i} $\gets$ \texttt{rowval[p]}
  \State \texttt{pos} $\gets$ \texttt{rowptr[i]}
  \State \texttt{colval[pos]} $\gets$ \texttt{j}
  \State \texttt{rowptr[i]} $\gets$ \texttt{rowptr[i] + 1}
  \EndFor
  \EndFor

  \For{$i$ $\gets$ \texttt{m} \textbf{down to} \texttt{2}}
  \State \texttt{rowptr[i]} $\gets$ \texttt{rowptr[i-1]}
  \EndFor
  \State \texttt{rowptr[1]} $\gets$ \texttt{1}
  \State \Return \texttt{rowptr}, \texttt{colval}
\end{algorithmic}

  \end{framed}
  \caption{Procedure to convert a sparse matrix $A$ stored in CSC format (\texttt{colptr}, \texttt{rowval}) into its equivalent CSR representation (\texttt{rowptr}, \texttt{colval}). The arrays use 1-based indexing.}
  \label{alg:transpose_sparsity_pattern}
\end{algorithm}

\begin{algorithm}[htp!]
  \footnotesize
  \begin{framed}
    \textbf{Input:}
\begin{itemize}
  \setlength{\itemsep}{0pt}
  \setlength{\parskip}{0pt}
  \setlength{\itemindent}{-12pt}
  \item \texttt{adjacency\_graph}, \texttt{vertices\_in\_order}, \texttt{nv}, \texttt{ne}
\end{itemize}

\textbf{Output:}
\begin{itemize}
  \setlength{\itemsep}{0pt}
  \setlength{\parskip}{0pt}
  \setlength{\itemindent}{-12pt}
  \item \texttt{color}, \texttt{star}, \texttt{hub}
\end{itemize}

\begin{algorithmic}[1]
  \State \texttt{color[1,\dots,nv]} $\gets$ \texttt{0}
  \State \texttt{forbidden\_colors[1,\dots,nv]} $\gets$ \texttt{0}
  \State \texttt{first\_neighbor[1,\dots,nv]} $\gets$ \texttt{(0,0,0)}
  \State \texttt{treated[1,\dots,nv]} $\gets$ \texttt{0}
  \State \texttt{star[1,\dots,ne]} $\gets$ \texttt{0}
  \State \texttt{hub} $\gets$ \texttt{[\,]}
  \State \texttt{num\_stars} $\gets$ \texttt{0}
  \For{each \texttt{v} $\in$ \texttt{vertices\_in\_order}}
  \For{each \texttt{(w, index\_vw}) $\in$ \texttt{indexed\_neighbors(adjacency\_graph, v)}}
  \If{\texttt{color[w]} $\ne$ \texttt{0}}
  \State \texttt{forbidden\_colors[color[w]]} $\gets$ \texttt{v}
  \State \texttt{(p, q, index\_pq)} $\gets$ \texttt{first\_neighbor[color[w]]}
  \If{\texttt{\texttt{p} \texttt{==} \texttt{v}}}
  \If{\texttt{treated[q]} $\neq$ \texttt{v}}
  \For{each \texttt{x $\in$ \texttt{neighbors(adjacency\_graph, q)}}}
  \If{\texttt{color[x]} $\ne$ \texttt{0}}
  \State \texttt{forbidden\_colors[color[x]]} $\gets$ \texttt{v}
  \EndIf
  \EndFor
  \State \texttt{treated[q]} $\gets$ \texttt{v}
  \EndIf
  \For{each \texttt{x $\in$ \texttt{neighbors(adjacency\_graph, w)}}}
  \If{\texttt{color[x]} $\ne$ \texttt{0}}
  \State \texttt{forbidden\_colors[color[x]]} $\gets$ \texttt{v}
  \EndIf
  \EndFor
  \State \texttt{treated[w]} $\gets$ \texttt{v}
  \Else
  \State \texttt{first\_neighbor[color[w]]} $\gets$ \texttt{(v, w, index\_{vw})}
  \For{each \texttt{(x, index\_wx}) $\in$ \texttt{indexed\_neighbors(adjacency\_graph, w)}}
  \If{\texttt{x} $\neq$ \texttt{v} \textbf{and} \texttt{color[x]} $\neq$ \texttt{0} \textbf{and} \texttt{x} \texttt{==} \texttt{hub[star[index\_{wx}]]}}
  \State \texttt{forbidden\_colors[color[x]]} $\gets$ \texttt{v}
  \EndIf
  \EndFor
  \EndIf
  \EndIf
  \EndFor
  \State \texttt{color[v]} $\gets$ \texttt{min\{i |} \texttt{forgiven\_colors[i]} $\ne$ \texttt{v}\}
  \For{each \texttt{(w, index\_vw}) $\in$ \texttt{indexed\_neighbors(adjacency\_graph, v)}}
  \If{\texttt{color[w]} $\ne$ \texttt{0}}
  \State \texttt{x\_exists} $\gets$ \texttt{false}
  \For{each \texttt{(x, index\_wx}) $\in$ \texttt{indexed\_neighbors(adjacency\_graph, w)}}
  \If{\texttt{x} $\neq$ \texttt{v} \textbf{and} \texttt{color[x] == color[v]}}
  \State \texttt{star\_{wx}} $\gets$ \texttt{star[index\_{wx}]}
  \State \texttt{hub[star\_{wx}]} $\gets$ \texttt{w}
  \State \texttt{star[index\_{vw}]} $\gets$ \texttt{star\_{wx}}
  \State \texttt{x\_exists} $\gets$ \texttt{true}
  \State \textbf{break}
  \EndIf
  \EndFor
  \If{(\textbf{not} \texttt{x\_exists})}
  \State \texttt{(p, q, index\_{pq})} $\gets$ \texttt{first\_neighbor[color[w]]}
  \If{\texttt{p} \texttt{==} \texttt{v} \textbf{and} \texttt{q} $\neq$ \texttt{w}}
  \State \texttt{star\_{vq}} $\gets$ \texttt{star[index\_{pq}]}
  \State \texttt{hub[star\_{vq}]} $\gets$ \texttt{v}
  \State \texttt{star[index\_{vw}]} $\gets$ \texttt{star\_{vq}}
  \Else
  \State \texttt{num\_stars} $\gets$ \texttt{num\_stars + 1}
  \State \texttt{push!(hub, -v)} \textbf{or} \texttt{push!(hub, -w)}
  \State \texttt{star[index\_{vw}]} $\gets$ \texttt{num\_stars}
  \EndIf
  \EndIf
  \EndIf
  \EndFor
  \EndFor
\end{algorithmic}

  \end{framed}
  \caption{Procedure \texttt{star\_coloring}.}
  \label{alg:star_coloring}
\end{algorithm}

\begin{algorithm}[ht!]
  \footnotesize
  \begin{framed}
    \textbf{Procedure} \texttt{create\_forest}(\texttt{ne})\\
$\bullet$ \textbf{Input:} An integer \texttt{ne} corresponding to the number of edges.\\
$\bullet$ \textbf{Output:} A structure \texttt{forest}.

\begin{algorithmic}[1]
  \item \texttt{nt} $\gets$ \texttt{ne}
  \item \texttt{parents[1..ne]} $\gets$ \texttt{[1, 2, ..., ne]}
  \item \texttt{ranks[1..ne]} $\gets$ \texttt{[0, 0, ..., 0]}
\end{algorithmic}

\medskip

\textbf{Procedure} \texttt{find\_root}(\texttt{forest}, \texttt{index})\\
$\bullet$ \textbf{Input:} One edge index, \texttt{index}, and the structure \texttt{forest}.\\
$\bullet$ \textbf{Output:} The root of the tree containing \texttt{index}.
\begin{algorithmic}[1]
  \If{\texttt{parents[index]} $\neq$ \texttt{index}}
  \State \texttt{parents[index]} $\gets$ \texttt{find\_root(parents[index], forest)}
  \EndIf
  \State \Return \texttt{parents[index]}
\end{algorithmic}

\medskip

\textbf{Procedure} \texttt{root\_union}(\texttt{forest}, \texttt{root1}, \texttt{root2})\\
$\bullet$ \textbf{Input:} Two root indices, \texttt{root1} and \texttt{root2}, and the structure \texttt{forest}.\\
$\bullet$ \textbf{Output:} The structure \texttt{forest} after merging the trees with roots \texttt{root1} and \texttt{root2}.
\begin{algorithmic}[1]
  \If{\texttt{ranks[root1]} \texttt{<} \texttt{ranks[root2]}}
  \State (\texttt{root1, root2}) $\gets$ (\texttt{root2, root1})
  \ElsIf{\texttt{ranks[root1]} \texttt{==} \texttt{ranks[root2]}}
  \State \texttt{ranks[root1]} $\gets$ \texttt{ranks[root1]} $+$ 1
  \EndIf
  \State \texttt{parents[root2]} $\gets$ \texttt{root1}
  \State \texttt{nt} $\gets$ \texttt{nt - 1}
\end{algorithmic}

  \end{framed}
  \caption{Routines related the structure \texttt{forest}.}
  \label{alg:forest}
\end{algorithm} 

\begin{algorithm}[ht!]
  \footnotesize
  \begin{framed}
    \textbf{Input:}
\begin{itemize}
  \setlength{\itemsep}{0pt}
  \setlength{\parskip}{0pt}
  \setlength{\itemindent}{-12pt}
  \item \texttt{adjacency\_graph}, \texttt{vertices\_in\_order}, \texttt{nv}, \texttt{ne}
\end{itemize}

\textbf{Output:}
\begin{itemize}
  \setlength{\itemsep}{0pt}
  \setlength{\parskip}{0pt}
  \setlength{\itemindent}{-12pt}
  \item \texttt{color}, \texttt{forest}
\end{itemize}

\begin{algorithmic}[1]
  \State \texttt{color[1,\dots,nv]} $\gets$ \texttt{0}
  \State \texttt{forbidden\_colors[1,\dots,nv]} $\gets$ \texttt{0}
  \State \texttt{first\_neighbor[1,\dots,nv]} $\gets$ \texttt{(0,0,0)}
  \State \texttt{first\_visit\_to\_tree[1,\dots,ne]} $\gets$ \texttt{(0,0)}
  \State \texttt{forest} $\gets$ \texttt{create\_forest(ne)}
  \For{each \texttt{v} $\in$ \texttt{vertices\_in\_order}}
  \For{each \texttt{(w, index\_vw}) $\in$ \texttt{indexed\_neighbors(adjacency\_graph, v)}}
  \If{\texttt{color[w]} $\neq$ 0}
  \State \texttt{forbidden\_colors[color[w]]} $\gets$ \texttt{v}
  \EndIf
  \EndFor
  \For{each \texttt{w $\in$ \texttt{neighbors(adjacency\_graph, v)}}}
  \If{\texttt{color[w]} $\neq$ 0}
  \For{each \texttt{(x, index\_wx}) $\in$ \texttt{indexed\_neighbors(adjacency\_graph, w)}}
  \If{\texttt{color[x]} $\neq$ 0 and \texttt{forbidden\_colors[color[x]]} $\neq$ \texttt{v}}
  \State \texttt{root\_wx} $\gets$ \texttt{find\_root(forest, index\_wx)}
  \State \texttt{(p, q)} $\gets$ \texttt{first\_visit\_to\_tree[root\_wx]}
  \If{\texttt{p} $\neq$ \texttt{v}}
  \State \texttt{first\_visit\_to\_tree[root\_wx]} $\gets$ \texttt{(v, w)}
  \ElsIf{\texttt{q} $\neq$ \texttt{w}}
  \State \texttt{forbidden\_colors[color[x]]} $\gets$ \texttt{v}
  \EndIf
  \EndIf
  \EndFor
  \EndIf
  \EndFor
  \State \texttt{color[v]} $\gets$ \texttt{min\{i |} \texttt{forgiven\_colors[i]} $\ne$ \texttt{v}\}
  \For{each \texttt{(w, index\_vw}) $\in$ \texttt{indexed\_neighbors(adjacency\_graph, v)}}
  \If{\texttt{color[w]} $\neq$ 0}
  \State \texttt{(p, q, index\_pq)} $\gets$ \texttt{first\_neighbor[color[w]]}
  \If{\texttt{p} $\neq$ \texttt{v}}
  \State \texttt{first\_neighbor[color[w]]} $\gets$ \texttt{(v, w, index\_vw)}
  \Else
  \State \texttt{root\_vw} $\gets$ \texttt{find\_root(forest, index\_vw)}
  \State \texttt{root\_pq} $\gets$ \texttt{find\_root(forest, index\_pq)}
  \State \texttt{root\_union(forest, root\_vw, root\_pq)}
  \EndIf
  \EndIf
  \EndFor
  \For{each \texttt{(w, index\_vw}) $\in$ \texttt{indexed\_neighbors(adjacency\_graph, v)}}
  \If{\texttt{color[x]} \texttt{==} \texttt{color[v]}}
  \State \texttt{root\_vw} $\gets$ \texttt{find\_root(forest, index\_vw)}
  \State \texttt{root\_wx} $\gets$ \texttt{find\_root(forest, index\_wx)}
  \If{\texttt{root\_vw} $\neq$ \texttt{root\_wx}}
  \State \texttt{root\_union(forest, root\_vw, root\_wx)}
  \EndIf
  \EndIf
  \EndFor
  \EndFor
\end{algorithmic}

  \end{framed}
  \caption{Procedure \texttt{acyclic\_coloring}.}
  \label{alg:acyclic_coloring}
\end{algorithm}

\begin{algorithm}[ht!]
  \footnotesize
  \begin{framed}
    \textbf{Input:}
\begin{itemize}
  \setlength{\itemsep}{0pt}
  \setlength{\parskip}{0pt}
  \setlength{\itemindent}{-12pt}
  \item \texttt{adjacency\_graph}, \texttt{nv}, \texttt{ne}, \texttt{nt}, \texttt{forest}
\end{itemize}

\textbf{Output:}
\begin{itemize}
  \setlength{\itemsep}{0pt}
  \setlength{\parskip}{0pt}
  \setlength{\itemindent}{-12pt}
  \item \texttt{nvmax}, \texttt{tree\_edge\_indices}, \texttt{tree\_vertices}, \texttt{tree\_neighbor\_indices}, \texttt{tree\_neighbors}
\end{itemize}

\begin{algorithmic}[1]
  \State \texttt{tree\_vertices[1,\dots,ne+nt]} $\gets$ \texttt{0}
  \State \texttt{tree\_neighbor\_indices[1,\dots,ne+nt+1]} $\gets$ \texttt{0}
  \State \texttt{tree\_neighbors[1,\dots,2*ne]} $\gets$ \texttt{0}
  \State \texttt{tree\_edge\_indices[1,\dots,nt+1]} $\gets$ \texttt{0}
  \State \texttt{root\_to\_tree[1,\dots,ne]} $\gets$ \texttt{0}
  \State \texttt{visited\_trees[1,\dots,nt]} $\gets$ \texttt{0}
  \State \texttt{vertex\_position[1,\dots,nt]} $\gets$ \texttt{0}
  \State \texttt{neighbor\_position[1,\dots,nt]} $\gets$ \texttt{0}
  \State \texttt{nr} $\gets$ 0
  \For{\texttt{index\_edge} $\gets$ 1 \textbf{to} \texttt{ne}}
  \State \texttt{root} $\gets$ \texttt{find\_root(forest, index\_edge)}
  \If{\texttt{root\_to\_tree[root]} \texttt{==} 0}
  \State \texttt{nr += 1}
  \State \texttt{root\_to\_tree[root]} $\gets$ \texttt{nr}
  \EndIf
  \State \texttt{index\_tree} $\gets$ \texttt{root\_to\_tree[root]}
  \State \texttt{tree\_edge\_indices[index\_tree + 1] += 1}
  \EndFor
  \State \texttt{nvmax} $\gets$ \texttt{maximum(tree\_edge\_indices)} $+$ 1
  \State \texttt{vertex\_position[1]} $\gets$ 0
  \State \texttt{neighbor\_position[1]} $\gets$ 0
  \For{\texttt{k} $\gets$ \texttt{2} \textbf{to} \texttt{nt}}
  \State \texttt{vertex\_position[k]} $\gets$ \texttt{vertex\_position[k-1] + tree\_edge\_indices[k] + 1}
  \State \texttt{neighbor\_position[k]} $\gets$ \texttt{neighbor\_position[k-1] + 2 * tree\_edge\_indices[k]}
  \EndFor
  \For{\texttt{j} $\gets$ \texttt{1} \textbf{to} \texttt{nv}}
  \For{each \texttt{(i, index\_ij}) $\in$ \texttt{indexed\_neighbors(adjacency\_graph, j)}}
  \If{\texttt{i} $\ne$ \texttt{j}}
  \State \texttt{root} $\gets$ \texttt{find\_root(forest, index\_ij)}
  \State \texttt{index\_tree} $\gets$ \texttt{root\_to\_tree[root]}
  \State \texttt{vertex\_index} $\gets$ \texttt{vertex\_position[index\_tree]}
  \If{\texttt{visited\_trees[index\_tree]} $\ne$ \texttt{j}}
  \State \texttt{visited\_trees[index\_tree]} $\gets$ \texttt{j}
  \State \texttt{vertex\_position[index\_tree] += 1}
  \State \texttt{vertex\_index += 1}
  \State \texttt{tree\_vertices[vertex\_index]} $\gets$ \texttt{j}
  \EndIf
  \State \texttt{neighbor\_position[index\_tree] += 1}
  \State \texttt{neighbor\_index} $\gets$ \texttt{neighbor\_position[index\_tree]}
  \State \texttt{tree\_neighbors[neighbor\_index]} $\gets$ \texttt{i}
  \State \texttt{tree\_neighbor\_indices[vertex\_index + 1] += 1}
  \EndIf
  \EndFor
  \EndFor
  \State \texttt{tree\_edge\_indices[1]} $\gets$ 1
  \State \texttt{tree\_neighbor\_indices[1]} $\gets$ 1
  \For{\texttt{k} = \texttt{2} \textbf{to} \texttt{nt+1}}
  \State \texttt{tree\_edge\_indices[k]} \!$\gets$\! \texttt{tree\_edge\_indices[k] \!+ tree\_edge\_indices[k-1]}
  \State \texttt{tree\_neighbor\_indices[k]} \!$\gets$\! \texttt{tree\_neighbor\_indices[k] \!+ tree\_neighbor\_indices[k-1]}
  \EndFor
\end{algorithmic}

  \end{framed}
  \caption{Procedure to recover the structure of the trees from \texttt{forest}.}
  \label{alg:structure_trees}
\end{algorithm}

\begin{algorithm}[ht!]
  \footnotesize
  \begin{framed}
    \textbf{Input:}
\begin{itemize}
  \setlength{\itemsep}{0pt}
  \setlength{\parskip}{0pt}
  \setlength{\itemindent}{-12pt}
  \item \texttt{nvmax} is the number of vertices in the largest tree of the forest.
  \item \texttt{tree\_edge\_indices} specifies the starting and ending indices of edges for each tree.
  \item \texttt{tree\_vertices} contains the list of vertices, grouped by tree.
  \item \texttt{tree\_neighbor\_indices} provides the positions of the first and last neighbors for each vertex in \texttt{tree\_vertices}, within the tree to which the vertex belongs.
  \item \texttt{tree\_neighbors} is a packed representation of the neighbors of each vertex in \texttt{tree\_vertices}.
\end{itemize}

\textbf{Output:}
\begin{itemize}
  \setlength{\itemsep}{0pt}
  \setlength{\parskip}{0pt}
  \setlength{\itemindent}{-12pt}
  \item \texttt{reverse\_bfs\_orders} is the edge ordering used for decompression, grouped by tree.
  \item \texttt{is\_star} is a boolean vector indicating whether each tree is a star.
\end{itemize}

\begin{algorithmic}[1]
  \State \texttt{reverse\_bfs\_orders[1,\dots,ne]} $\gets$ \texttt{(0,0)}
  \State \texttt{degrees[1,\dots,nv]} $\gets$ \texttt{0}
  \State \texttt{reverse\_mapping[1,\dots,nv]} $\gets$ \texttt{0}
  \State \texttt{queue[1,\dots,nvmax]} $\gets$ \texttt{0}
  \State \texttt{is\_star[1,\dots,nt]} $\gets$ \texttt{false}
  \State \texttt{num\_edges\_treated} $\gets$ 0
  \For{$k = 1$ \textbf{to} \texttt{nt}}
  \State \texttt{queue\_start} $\gets$ 1
  \State \texttt{queue\_end} $\gets$ 0
  \State \texttt{first\_vertex} $\gets$ \texttt{tree\_edge\_indices[k] + k - 1}
  \State \texttt{last\_vertex} $\gets$ \texttt{tree\_edge\_indices[k + 1] + k - 1}
  \For{\texttt{index\_vertex} \textbf{in} \texttt{first\_vertex:last\_vertex}}
  \State \texttt{vertex} $\gets$ \texttt{tree\_vertices[index\_vertex]}
  \State \texttt{index\_first\_neighbor} =  \texttt{tree\_neighbor\_indices[index\_vertex]}
  \State \texttt{degree} = \texttt{tree\_neighbor\_indices[index\_vertex+1]} - \texttt{index\_first\_neighbor}
  \State \texttt{degrees[vertex]} $\gets$ \texttt{degree}
  \State \texttt{reverse\_mapping[vertex]} $\gets$ \texttt{index\_vertex}
  \If{\texttt{degree} \texttt{==} 1}
  \State \texttt{queue\_end += 1}
  \State \texttt{queue[queue\_end]} $\gets$ \texttt{vertex}
  \EndIf
  \EndFor
  \State \texttt{nv\_tree} $\gets$ \texttt{tree\_edge\_indices[k + 1] - tree\_edge\_indices[k] + 1}
  \State \texttt{is\_star[k]} $\gets$ \texttt{queue\_end} $\ge$ \texttt{nv\_tree - 1}
  \While{\texttt{queue\_start} $\le$ \texttt{queue\_end}}
  \State \texttt{leaf} $\gets$ \texttt{queue[queue\_start]}
  \State \texttt{queue\_start += 1}
  \State \texttt{degrees[leaf]} $\gets$ 0
  \State \texttt{index\_leaf} $\gets$ \texttt{reverse\_mapping[leaf]}
  \State \texttt{first\_neighbor} $\gets$ \texttt{tree\_neighbor\_indices[index\_leaf]}
  \State \texttt{last\_neighbor} $\gets$ \texttt{tree\_neighbor\_indices[index\_leaf + 1] - 1}
  \For{\texttt{index\_neighbor} \textbf{in} \texttt{first\_neighbor:last\_neighbor}}
  \State \texttt{neighbor} $\gets$ \texttt{tree\_neighbors[index\_neighbor]}
  \If{\texttt{degrees[neighbor]} $\ne$ 0}
  \State \texttt{num\_edges\_treated += 1}
  \State \texttt{reverse\_bfs\_orders[num\_edges\_treated]} $\gets$ (\texttt{leaf}, \texttt{neighbor})
  \State \texttt{degrees[neighbor] -= 1}
  \If{\texttt{degrees[neighbor]} \texttt{==} 1}
  \State \texttt{queue\_end += 1}
  \State \texttt{queue[queue\_end]} $\gets$ \texttt{neighbor}
  \EndIf
  \EndIf
  \EndFor
  \EndWhile
  \EndFor
\end{algorithmic}

  \end{framed}
  \caption{Procedure the compute the reverse BFS order of the trees.}
  \label{alg:reverse_bfs}
\end{algorithm}

\begin{algorithm}[ht!]
  \footnotesize
  \begin{framed}
    \textbf{Input:}
\begin{itemize}
  \setlength{\itemsep}{0pt}
  \setlength{\parskip}{0pt}
  \setlength{\itemindent}{-12pt}
  \item \texttt{buckets}: a vector of stacks, one per degree
  \item \texttt{v}: a vertex to be moved up one degree
\end{itemize}

\begin{algorithmic}[1]
  \State \texttt{d} $\gets$ current dynamic degree of \texttt{v}
  \State \texttt{p} $\gets$ position of \texttt{v} inside \texttt{buckets[d]}
  \State \texttt{w} $\gets$ \texttt{buckets[d][end]}
  \State \texttt{buckets[d][p]} $\gets$ \texttt{w}
  \State \texttt{buckets[d][end]} $\gets$ \texttt{v}
  \State \texttt{pop\_last(buckets[d])}
  \State \texttt{append\_last(buckets[d+1], v)}
  \State record new positions of \texttt{v} and \texttt{w}
\end{algorithmic}

  \end{framed}
  \caption{Stack-based degree bucket update, following \texttt{ColPack}.}
  \label{alg:increase_degree_stacks}
\end{algorithm}

\begin{algorithm}[ht!]
  \footnotesize
  \begin{framed}
    \textbf{Input:}
\begin{itemize}
  \setlength{\itemsep}{0pt}
  \setlength{\parskip}{0pt}
  \setlength{\itemindent}{-12pt}
  \item \texttt{bucketvec}: a vector of integers, one per vertex
  \item \texttt{bucketlims}: a vector of integer pairs, one per degree
  \item \texttt{v}: a vertex to be moved up one degree
\end{itemize}

\begin{algorithmic}[1]
  \State \texttt{d} $\gets$ current dynamic degree of \texttt{v}
  \State \texttt{p} $\gets$ position of \texttt{v} inside \texttt{buckets[d]}
  \State \texttt{(a, b)} $\gets$ \texttt{bucketlims[d]}
  \State \texttt{(a\_up, b\_up)} $\gets$ \texttt{bucketlims[d+1]}
  \State \texttt{w} $\gets$ \texttt{buckets[d][end]}
  \State \texttt{bucketvec[a+p]} $\gets$ \texttt{w}
  \State \texttt{bucketvec[a\_up-1]} $\gets$ \texttt{v}
  \State \texttt{bucketlims[d]} $\gets$ \texttt{(a, b-1)}
  \State \texttt{bucketlims[d+1]} $\gets$ \texttt{(a\_up-1, b\_up)}
  \State record new positions of \texttt{v} and \texttt{w}
\end{algorithmic}

  \end{framed}
  \caption{Vector-based degree bucket update, introduced in \texttt{SMC}.}
  \label{alg:increase_degree_vector}
\end{algorithm}


\end{document}